\documentclass[11pt]{amsart} 
\usepackage{amssymb}
\usepackage{amsmath}
\usepackage{mathtools}
\usepackage{enumitem}
\usepackage{mathrsfs}
\usepackage[english]{babel}
\usepackage[all]{xy}
\usepackage{hyperref}
\usepackage{marginnote}
\usepackage{comment}
\usepackage{tikz-cd}

\usepackage{stmaryrd}

\usepackage{fullpage}

\RequirePackage{mathrsfs} 

\newtheorem{theorem}{Theorem}[section]
\newtheorem{lemma}[theorem]{Lemma}
\newtheorem{proposition}[theorem]{Proposition}
\newtheorem{corollary}[theorem]{Corollary}
\newtheorem{conjecture}[theorem]{Conjecture}
\newtheorem{problem}[theorem]{Problem}
\newtheorem{alphtheorem}{Theorem}
\newtheorem{alphremark}[alphtheorem]{Remark}

\theoremstyle{definition}
\newtheorem*{ack}{Acknowledgments}

\newtheorem{remark}[theorem]{Remark}

\newtheorem{definition}[theorem]{Definition}

\numberwithin{equation}{section} \numberwithin{figure}{section}

\DeclareMathOperator{\Pic}{Pic} 
\DeclareMathOperator{\Gal}{Gal} 
\DeclareMathOperator{\Aut}{Aut} 

\DeclareMathOperator{\Sym}{Sym}
\DeclareMathOperator{\Spec}{Spec}

\DeclareMathOperator{\an}{an}

\DeclareMathOperator{\sat}{sat}
\DeclareMathOperator{\Hilb}{Hilb}

\DeclareMathOperator{\id}{id}
\DeclarePairedDelimiter\abs{\lvert}{\rvert}

\newcommand{\Qbar}{\overline{\QQ}}

\newcommand\ZZ{\mathbb{Z}}

\newcommand\QQ{\mathbb{Q}}

\newcommand\CC{\mathbb{C}}

\newcommand*\ratmap{\mathbin{\tikz [baseline=0ex,-latex, dashed, ->] \draw [densely dashed] (0em,0.58ex) -- (1.3em,0.58ex);}}

\usepackage{color}

\definecolor{orange}{rgb}{1,0.5,0}

\title[Symmetric products and puncturing Campana-special varieties]{Symmetric products and puncturing Campana-special varieties}

\author{Finn Bartsch}
\address{Finn Bartsch \\
IMAPP Radboud University Nijmegen \\
PO Box 9010, 6500GL \\
Nijmegen, The Netherlands\\}
\email{f.bartsch@math.ru.nl }
 
\author{Ariyan Javanpeykar}
\address{Ariyan Javanpeykar \\ 
IMAPP Radboud University Nijmegen \\
PO Box 9010, 6500GL\\
 Nijmegen, The Netherlands}
\email{ariyan.javanpeykar@ru.nl }
 
\author{Aaron Levin}
\address{Aaron Levin\\
Department of Mathematics \\
Michigan State University \\
619 Red Cedar Road\\
East Lansing, MI 48824\\
USA.}
\email{levina@msu.edu}



\begin{document}

\begin{abstract} 
We give a counterexample to the Arithmetic Puncturing Conjecture and Geometric Puncturing Conjecture of Hassett--Tschinkel using symmetric powers of uniruled surfaces, and propose a corrected conjecture inspired by Campana's conjectures on special varieties. We confirm Campana's conjecture on potential density for symmetric powers of products of curves. As a by-product, we obtain an example of a surface without a potentially dense set of rational points, but for which some symmetric power does have a dense set of rational points, and even satisfies Corvaja--Zannier's version of the Hilbert property.
\end{abstract}

\maketitle
\thispagestyle{empty}

\section{Introduction}
The aim of this paper is to give a counterexample to the Puncturing Conjectures of Hassett--Tschinkel using symmetric products of surfaces, and to propose corrected conjectures guided by Campana's conjectures on special varieties, dense entire curves, and potential density of rational points over number fields and function fields, respectively. 

We start with an overview of Campana's conjectures for quasi-projective varieties. To do so, let $k$ be an algebraically closed field of characteristic zero. A variety over $k$ is a finite type separated integral scheme over $k$. 

Central to this paper is the class of special varieties introduced by Campana in \cite{Ca04} for smooth projective varieties and \cite[Definition~8.1]{Ca11} in his more general orbifold setting. 
 We give a definition here, and refer to Section \ref{section:puncturing_specials} for a discussion of some basic properties of special varieties.
A pair $(X,D)$ is an \emph{snc pair} if $X$ is a smooth proper variety over $k$ and $D$ is a simple normal crossings divisor on $X$. 
We follow \cite[\S 11]{Iitaka} and let $\Omega^1_X(\log D)\subset \Omega^1_X$ be the subsheaf of differential forms with log poles along $D$. Define $\Omega^p_X(\log D) = \Lambda^p \Omega^1_X(\log D)$. Bogomolov showed that for every line bundle $\mathcal{L}$ admitting a nonzero morphism $\mathcal{L}\to \Omega^p_X(\log D)$, the Iitaka dimension $\kappa(\mathcal{L})$ is at most $p$; see \cite[\S 12, Theorem 4]{BogomolovHol} for the projective case and \cite[Corollary~6.9]{EVbook} in the snc case. For snc pairs, the following definition encapsulates all we need.

\begin{definition}\label{def:special_pair}
Let $(X,D)$ be an snc pair. For $1 \leq p \leq \dim X$, a line bundle $\mathcal{L}$ on $X$ is a \emph{Bogomolov sheaf of rank $p$ (for $(X,D)$)} if there is a nonzero morphism $\mathcal{L} \to \Omega^p_X(\log D)$ and the Iitaka dimension $\kappa(\mathcal{L})$ of $\mathcal{L}$ is equal to $p$. A line bundle $\mathcal{L}$ on $X$ is a \emph{Bogomolov sheaf (for $(X,D)$)} if there is an integer $1 \leq p \leq \dim X$ such that $\mathcal{L}$ is a Bogomolov sheaf of rank $p$. The snc pair $(X,D)$ is \emph{special} if it does not have any Bogomolov sheaves. 
\end{definition}

For (possibly very singular) varieties, the notion of specialness is defined by passing to an snc model. More precisely:

\begin{definition}
A variety $X$ over $k$ is \emph{special} if there is a resolution of singularities $Y \to X$ and a smooth projective compactification $\overline{Y}$ of $Y$ whose boundary $\overline{Y}\setminus Y =: D$ is an snc divisor such that the snc pair $(\overline{Y},D)$ is special.
\end{definition}

By Lemma \ref{lemma:independence_of_choices}, this definition is independent of the choice of the resolution and compactification.

\subsection{Complex-analytic notions of specialness}
We now introduce the conjecturally equivalent counterparts to Campana's notion of specialness.
 
\begin{definition}[Brody specialness]
A variety $X$ over $\mathbb{C}$ is \emph{Brody-special} if there is a holomorphic map $\mathbb{C}\to X^{\an}$ whose image is Zariski-dense in $X$. 
\end{definition}
 
If $X$ is a variety over $\mathbb{C}$, we let $d_X$ denote the Kobayashi pseudometric on $X^{\an}$. 
This pseudometric plays a crucial role in Campana's conjecture through the following notion (which Campana refers to as ``hyperbolically special'' \cite[Definition~9.1.1]{Ca04}).

\begin{definition}[Kobayashi-specialness]
A variety $X$ over $\mathbb{C}$ is \emph{Kobayashi-special} if there is a proper birational morphism $X'\to X$ such that $X'$ is a smooth variety with $d_{X'} \equiv 0$. 
\end{definition}

It follows from a classical theorem of Campbell--Ogawa and Campbell--Howard--Ochiai (see Theorem \ref{thm:punc_koba} below) that, if $X$ is Kobayashi-special, then the Kobayashi pseudometric $d_Y$ is identically zero for any desingularization $Y\to X$.

By the distance-decreasing property of the Kobayashi pseudometric, if $X$ is a Kobayashi-special variety over $\mathbb{C}$, then $d_X\equiv 0$ (as any desingularization $X'\to X$ is surjective and has vanishing pseudometric). However, the condition that $d_X \equiv 0$ does not necessarily imply that $X$ is Kobayashi-special if $X$ is singular. For example, the cone over a hyperbolic curve has vanishing Kobayashi pseudometric (as it is covered by different copies of $\mathbb{G}_m$), but it is not Kobayashi-special as it dominates (up to blow-up) a hyperbolic curve. This shows that the notion of Kobayashi-specialness really requires passing to a desingularization. 

Campana conjectured that the above three notions are all equivalent:
\begin{conjecture}[Campana] \label{conj:campana_analytic}
Let $X$ be a variety over $\mathbb{C}$. Then the following are equivalent.
\begin{enumerate}
\item $X$ is special.
\item $X$ is Brody-special.
\item $X$ is Kobayashi-special.
\end{enumerate}
\end{conjecture}

Although this conjecture is stated for all varieties, it easily reduces to the smooth case.

\subsection{Arithmetic specialness}
Recall that $k$ is an algebraically closed field of characteristic zero.

The arithmetic property that should characterize the property of being special for a variety $X$ over $k$ is that there is an abundance of rational points on $X$. To make this more precise, let $S$ be an integral noetherian scheme with function field $K:=K(S)$ and let $X \to S$ be a morphism of schemes. We define $X(S)^{(1)}$ to be the set of $P$ in $X(K)$ such that, for every point $s$ in $S$ of codimension one, the point $P$ lies in the image of $X(\mathcal{O}_{S,s}) \to X(K)$.
Vojta refers to the points in $X(S)^{(1)}$ as \emph{near-integral $S$-points}; see \cite{VojtaExc}. If $S$ is one-dimensional, then $X(S) = X(S)^{(1)}$, so that near-integral $S$-points are the same as $S$-points. Moreover, if $X\to S$ is proper, then $X(S)^{(1)} = X(K) = X_K(K)$, i.e., the $K$-rational points on $X_K$ are the near-integral $S$-points of $X$. The notion of near-integral $S$-points is the ``correct'' notion to consider when studying rational points on proper varieties over finitely generated fields of positive transcendence degree over $\mathbb{Q}$.

\begin{definition}[Arithmetic specialness]
A variety $X$ over $k$ is \emph{arithmetically-special over $k$} if there is a $\ZZ$-finitely generated subring $A\subset k$ and a finite type separated model $\mathcal{X}$ for $X$ over $A$ such that the set $\mathcal{X}(A)^{(1)}$ of near-integral $A$-points is dense in $X$.
\end{definition}
 
For example, a variety $X$ over $\Qbar$ is arithmetically-special over $\Qbar$ if and only if there is a number field $K$, a finite set of finite places $S$ of $K$, and a model $\mathcal{X}$ for $X$ over $\mathcal{O}_{K,S}$ such that $\mathcal{X}(\mathcal{O}_{K,S})$ is dense in $X$. Moreover, a proper variety $X$ over $\Qbar$ (resp. $k$) is arithmetically-special over $\Qbar$ (resp. $k$) if and only if there is a number field $K\subset \Qbar$ (resp. a finitely generated subfield $K\subset k$) and a proper model $\mathcal{X}$ for $X$ over $K$ such that $\mathcal{X}(K)$ is dense in $X$. 
 
Arithmetic specialness is a formal way of capturing the well-studied property of having a potentially dense set of rational points.
Examples of arithmetically-special varieties include curves of genus at most one, unirational varieties, abelian varieties, Enriques surfaces, certain K3 surfaces (and conjecturally all), and certain Fano varieties (and, again, conjecturally all) \cite{HassettSurvey, HassettTschinkel}. 

One of our main results is that certain symmetric products of non-arithmetically-special surfaces are arithmetically-special (see Theorem \ref{thm:sym_prods_1}).

\subsection{Geometric specialness}
A function field analogue of the notion of arithmetic specialness was introduced in \cite{JR}. Roughly speaking, instead of asking for the abundance of rational points over a number field, one asks for the abundance of pointed curves (which figure as rational points over function fields).

\begin{definition}[Geometrically-special]\label{def:geom_specialness}
A variety $X$ over $k$ is \emph{geometrically-special over $k$} if, for every dense open subset $U\subset X$, there exists a smooth affine connected pointed curve $(C,c)$, a point $x$ in $U(k)$, and a sequence of morphisms $\{f_i \colon (C,c)\to (X,x)\}_{i=1}^\infty$ such that $C \times X$ is covered by the graphs $\Gamma_{f_i}$ of these maps, i.e., the closure of $\cup_{i=1}^{\infty} \Gamma_{f_i}$ in $C \times X$ equals $C \times X$. 
\end{definition}

For a variety to be geometrically-special means, roughly speaking, that it is covered by curves in a particularly strong sense. This notion was studied (mostly for projective varieties) in \cite{JR}, but also \cite{BJR, PRT}. The following conjecture is essentially due to Campana.
 
\begin{conjecture}[Campana] \label{conj:arithmetic}
Let $X$ be a variety over an algebraically closed field $k$ of characteristic zero. Then the following are equivalent.
\begin{enumerate}
\item $X$ is special.
\item $X$ is arithmetically-special over $k$.
\item $X$ is geometrically-special over $k$.
\end{enumerate}
\end{conjecture}

None of the implications are known in full generality, unless $X$ is one-dimensional or a closed subvariety of an abelian variety. Indeed, if $X$ is a closed subvariety of an abelian variety, then the above conjecture follows from the work of Faltings and Yamanoi \cite{FaltingsLang, Yamanoi1} (see \cite[Theorem~3.5]{JR} for a detailed explanation).

Conjectures \ref{conj:campana_analytic} and \ref{conj:arithmetic} provide a plethora of predictions, and the aim of this paper is to investigate predictions made for smooth varieties deprived of a closed subset of codimension at least two. There are other aspects of Campana's conjectures pertaining to non-archimedean specialness \cite{MorrowRosso} and numerical dimension \cite{PRT, Wu} which we do not discuss here. 

Guided by these predictions, we prove that certain symmetric products of non-geometrically-special surfaces are geometrically-special (see Theorem \ref{thm:sym_prods_1}).

\subsection{Hilbert irreducibility}

Campana's arithmetic conjecture predicts that a special variety has a potentially dense set of integral points over some suitable $\mathbb{Z}$-finitely generated subring of $k$. In other words, special varieties should have many integral points. Quantifying what ``many'' points could mean (besides mere density) naturally leads to Hilbert-type properties (studied originally for their relation to the inverse Galois problem \cite[\S~3]{SerreTopicsGalois}).

We follow \cite{CZHP} (see also \cite[Definition~1.2]{CDJLZ}) and introduce the weak Hilbert property. Note that a morphism $Y\to X$ of normal (geometrically integral) varieties is a \emph{ramified cover} if it is finite surjective and not \'etale. 

\begin{definition}[Corvaja--Zannier] \label{def:whp}
A normal proper variety $X$ over a field $K$ has the \emph{weak Hilbert property over $K$} if for every integer $n\geq 1$ and every finite collection $(\pi_i \colon Y_i \to X)$ of ramified covers with each $Y_i$ a normal variety over $K$, the set 
\[ X(K) \setminus \cup_{i=1}^n \pi_i (Y_i(K)) \]
is dense in $X$.
\end{definition}

Note that in our definition of the weak Hilbert property we consider proper varieties (hence $K$-points) only for simplicity's sake; the more general definition for quasi-projective schemes over regular $\mathbb{Z}$-finitely generated subrings of $k$ is given by Luger in \cite[Definition~1.3]{Luger2}. 

In the study of liftabilty of rational points along ramified covers of not necessarily proper varieties (e.g., punctured varieties), the notion of a strongly thin subset is indispensable: 

\begin{definition} 
Let $X$ be a normal variety over a field $K$. A subset $\Omega\subset X(K)$ is \emph{strongly thin} if there is an integer $n\geq 1$ and a finite collection $(\pi_i \colon Y_i \to X)$ of finite ramified covers such that 
\[ \Omega \setminus \cup_{i=1}^n \pi_i (Y_i(K)) \]
is not dense.
\end{definition}

With this definition, a normal proper variety $X$ over a field $k$ has the weak Hilbert property over $k$ (in the sense of Definition \ref{def:whp}) if and only if $X(k)$ is not strongly thin.

In the non-proper setting, the definition of the weak Hilbert property pertains to density of near-integral points (as in the definition of an arithmetically-special variety).

\begin{definition}
Let $k$ be an algebraically closed field of characteristic zero. A normal variety $X$ over $k$ has the \emph{arithmetic weak Hilbert property over $k$} if there is a $\ZZ$-finitely generated subring $A\subset k$ and a finite type separated model $\mathcal{X}$ for $X$ over $A$ such that the set $\mathcal{X}(A)^{(1)}$ of near-integral $A$-points is not strongly thin in $X$.
\end{definition} 

Given a proper variety $X$ over a number field $K$, the variety $X_{\overline{K}}$ has the arithmetic weak Hilbert property over $\overline{K}$ if, for every finite collection $(\pi_i \colon Y_i \to X)$ of finite ramified covers, the set $X(K)\setminus \cup_{i=1}^n \pi_i (Y_i(K))$ is dense.

Note that if $X$ has the arithmetic weak Hilbert property over $k$ and $L/k$ is an extension of algebraically closed fields, then $X_L$ has the arithmetic weak Hilbert property over $L$; this is a consequence of \cite[Proposition~3.2]{BSFPextensions}. Obviously, if $X$ has the arithmetic weak Hilbert property, then $X$ is arithmetically-special. 

The weak Hilbert property for $X$ means, roughly speaking, that any ramified cover of $X$ has ``fewer'' points than $X$ (and that $X$ has many points itself). In \cite{CZHP} Corvaja--Zannier conjectured that any smooth projective variety with a potentially dense set of rational points has the weak Hilbert property potentially (we note that the smoothness assumption here is crucial, see Remark \ref{remark:enriques}). Combined with Campana's conjecture (Conjecture \ref{conj:arithmetic}) in the general quasi-projective setting this leads to the following:

\begin{conjecture}[Campana, Corvaja--Zannier] \label{conj:whp}
Let $X$ be a smooth variety over a finitely generated field $K$ of characteristic zero. Then the following are equivalent.
\begin{enumerate}
\item The variety $X_{\overline{K}}$ is special.
\item The variety $X_{\overline{K}}$ is arithmetically-special. 
\item The variety $X_{\overline{K}}$ has the arithmetic weak Hilbert property. 
\end{enumerate}
\end{conjecture}

The weak Hilbert property also has ``Brody'' (resp. ``Kobayashi'', resp. ``geometric'') analogues which are probably equivalent to Brody-specialness (resp. Kobayashi-specialness, resp. geometric specialness); see, for example, \cite{CampanaWinkelmann} and \cite{CampanaBook-2}. We omit a further discussion of these topics here, and focus primarily on the arithmetic aspects.

\subsection{Symmetric products} 
If $X$ is a quasi-projective variety and $m\geq 1$ is an integer, the permutation group $S_m$ acts on $X^m$ by $\sigma\cdot (x_1,\ldots,x_m) = (x_{\sigma(1)},\ldots, x_{\sigma(m)})$. We let $\Sym^m(X) = X^m/S_m$ denote the $m$-th symmetric power of $X$.

Note that $\Sym^m(X)$ is an $m\cdot \dim(X)$-dimensional quasi-projective variety, and that the quotient morphism $X^m\to\Sym^m(X)$ is a finite surjective morphism. Let $\Delta_{i,j}\subset X^m$ be the closed subscheme given by the set of points $P = (x_1,\ldots,x_m)\in X^m$ satisfying $x_i = x_j$. We define the \emph{big diagonal} $\Delta^{(m)}$ of the $m$-th symmetric power $\Sym^m(X)$ to be the image of $\Delta^m:=\cup_{1\leq i < j \leq m} \Delta_{i,j}$. If $X$ is smooth, the closed subset $\Delta^{(m)}$ contains the singular locus of $\Sym^m(X)$, as the morphism $X^m\setminus \Delta^m \to \Sym^m(X)\setminus \Delta^{(m)}$ is an $S_m$-torsor, and thus finite \'etale. 

Arapura and Archava \cite{ArapuraArchava} showed that any symmetric power of a general type variety of dimension at least two is of general type. Conversely, if the symmetric power of a variety is of general type, then obviously the variety itself is of general type. It is natural to ask whether similar statements hold for the antithesis of the class of varieties of general type (i.e., the class of special varieties). It is not hard to show that, if $X$ is special, all of its symmetric powers will be special. However, it can very well happen that the symmetric power of a non-special variety is special. Let us be more precise.

Let $C$ be a smooth projective connected curve of genus $g\geq 2$ over $k$, and let $m\geq g$ be an integer. 
Central to this paper are the (singular!) symmetric powers of the surface $C\times \mathbb{P}^1_k$. As shown in \cite[Theorem~3]{CCR}, we have the following result pertaining to Campana's Conjecture \ref{conj:campana_analytic}.

\begin{theorem}[Campana--Cadorel--Rousseau]\label{thm:sym_prods_ccr}
Let $C$ be a smooth projective connected curve of genus $g$ over $\mathbb{C}$, and let $m\geq g$ be a positive integer. Then the following statements hold.
\begin{enumerate}
\item The variety $\Sym^m(C\times \mathbb{P}^1)$ is special.
\item The variety $\Sym^m(C\times \mathbb{P}^1)$ is Brody-special. 
\item The variety $\Sym^m(C\times \mathbb{P}^1)$ is Kobayashi-special.
\end{enumerate}
\end{theorem}

Guided by Conjecture \ref{conj:arithmetic} we verify that the special variety $\Sym^m(C\times \mathbb{P}^1)$ is both arithmetically-special and geometrically-special (see Corollaries \ref{cor:sym_n_is_geom_special} and \ref{cor:sym_n_is_arith_special} below).

\begin{alphtheorem}\label{thm:sym_prods_1}
Let $C$ be a smooth projective connected curve of genus $g$ over $k$, and let $m\geq g$ be a positive integer. Then $\Sym^m(C\times \mathbb{P}^1)$ is arithmetically-special over $k$ and geometrically-special over $k$.
\end{alphtheorem}

The study of potential density of rational points on symmetric powers $\Sym^n(X)$ of a surface $X$ is not new. 
For example, in \cite{HassettTschinkel}, it is shown that the Kodaira dimension of $\Sym^n(X)$ is $n$ times the Kodaira dimension of $X$. This leads Hassett and Tschinkel to predict that the behaviour of rational points on $X$ and $\Sym^n(X)$ should be similar (see \cite[p.~2]{HassettTschinkel}). Note that Theorem \ref{thm:sym_prods_1} contradicts this expectation.

Motivated by Corvaja--Zannier's conjectures on the Hilbert property, we also establish the stronger fact that $\Sym^m(C\times \mathbb{P}^1)$ has the potential weak Hilbert property (see Theorem \ref{thm:whp_for_sym} below).

\begin{alphtheorem}\label{thm:whp_for_sym_intro} 
Let $C$ be a smooth projective curve of genus $g$ over a finitely generated field $K$ of characteristic zero and $m \geq g$ a positive integer. Then there is a finite field extension $L/K$ such that $\Sym^m(C_L \times \mathbb{P}^1_L)$ has the weak Hilbert property over $L$.
\end{alphtheorem}

Our proof of Theorem \ref{thm:whp_for_sym_intro} uses the recently established version of Hilbert's irreducibility theorem for abelian varieties \cite{CDJLZ}. In fact, to prove the (potential) weak Hilbert property for $\Sym^m(C\times \mathbb{P}^1_K)$, we first establish a version of Hilbert's irreducibility theorem for the symmetric product $\Sym^m(C)$ of the curve $C$; this leads to an interesting application pertaining to the infinitude of $S_m$-Galois points on $C$ (see Corollary \ref{cor:Sm_points_on_C} for a precise statement).

Note that for $E$ an elliptic curve and $C$ as in Theorem \ref{thm:whp_for_sym}, the variety $\Sym^m(C\times E)$ is special and Brody-special \cite[Theorem~3]{CCR}. However, we are surprisingly not able to prove that it has a dense set of rational points over a large enough number field, unless $C$ dominates $E$. The situation is similar in the (isotrivial) function field setting: we are only able to prove that $\Sym^m(C\times E)$ is geometrically-special if $C$ dominates $E$ (see Theorems \ref{thm:symCcrossE_geom_special} and \ref{thm:symCcrossE_arith_special} below).

\begin{alphtheorem}\label{thm:symCcrossE}
Let $C$ be a smooth projective curve of genus $g$ over a finitely generated field $K$ of characteristic zero. Let $E$ be an elliptic curve over $K$ such that $C_{\overline{K}}$ dominates $E_{\overline{K}}$. If $m \geq g$, then $\Sym^m(C_{\overline{K}} \times E_{\overline{K}})$ is arithmetically-special and geometrically-special. 
\end{alphtheorem}

In the proof of Theorem \ref{thm:symCcrossE}, we invoke the following criterion for density of graphs which is established using properties of Hilbert schemes (see Theorem \ref{thm:geomspeccriterion}); we believe this density criterion to be of independent interest.

\begin{alphtheorem}\label{thm:geomspeccriterion_intro}  
Let $Y$ be a variety over $k$ and let $X$ be a quasi-projective variety. Let $(\phi_i \colon Y \to X)_{i \in I}$ be a family of morphisms. Suppose that there is a point $y_0 \in Y$ such that $\{ \phi_i(y_0)~|~i \in I\}$ is dense in $X$. Then the union of graphs $\bigcup \Gamma_{\phi_i}$ is dense in $Y \times X$.
\end{alphtheorem} 

It remains an open problem to show that $\Sym^m(C\times E)$ is arithmetically-special (resp. geo\-me\-tri\-cal\-ly-special), even for $g=m=2$. If $g=m=2$, enlarging the base field $K$ appropriately, we are naturally led to investigate whether there is a collection of quadratic points $c_1, c_2,\ldots \in C$ such that the associated collection in $\Sym^2(C)(K)$ is dense and such that, for every $i=1,2,\ldots$, the rank of $E$ over the residue field $K(c_i)$ of $c_i$ is strictly larger than the rank of $E(K)$. However, we do not know how to prove the existence of such a collection of quadratic points. 

Note that in this paper we are mostly concerned with the specialness of symmetric powers of non-special varieties. It is however also natural to study the hyperbolicity of such symmetric powers. For example, if $S$ is a smooth projective hyperbolic variety  over $\mathbb{C}$, then one can show that $\Sym^m(S)$  is also hyperbolic, under suitable assumptions (see \cite{CCR, GFP}). 
 
If $X$ is a special (resp. arithmetically-special) variety over $k$, then it is obvious that $\Sym^m(X)$ is special (resp. arithmetically-special). Indeed, in the arithmetic setting, if $X$ has a dense set of integral points, then so does $X^m$. Projecting these integral points along $X^m \to \Sym^m(X)$, it follows directly that $\Sym^m(X)$ has a dense set of integral points as well. On the other hand, it is not at all clear whether some smooth model of $\Sym^m(X)$ has the arithmetic weak Hilbert property; note that Conjecture \ref{conj:whp} predicts that this is the case. If $X$ is rational over $K$, then $\Sym^m(X)$ is rational as well \cite{Mattuck} and thus satisfies the Hilbert property. Moreover, if $X$ is a smooth projective rationally connected variety satisfying a certain strong form of weak approximation, then $\Sym^m(X)$ does as well \cite[Theorem~1.3]{ChenZhang}; in particular, for such $X$, some smooth model of $\Sym^m(X)$ has the Hilbert property. However, we do not know whether some smooth model of $\Sym^m(A)$ satisfies the potential weak Hilbert property if $A$ is an abelian variety of dimension at least two.

\subsection{The Puncturing Problems}
In Problem 2.11 and Problem 2.14 of \cite{HT}, Hassett and Tschinkel proposed the following ``Arithmetic Puncturing Problem'' and ``Geometric Puncturing Problem'':

\begin{problem}[Arithmetic Puncturing Problem] \label{pAP}
Let $X$ be a projective variety with canonical singularities and $Z$ a subvariety of codimension $\geq 2$. Assume that rational points on $X$ are potentially dense. Are integral points on $(X,Z)$ potentially dense? (In other words, if $X$ is arithmetically-special, is $X\setminus Z$ also arithmetically-special?)
\end{problem}

\begin{problem}[Geometric Puncturing Problem] \label{pGP}
Let $X$ be a projective variety with canonical singularities and $Z$ a subvariety of codimension $\geq 2$. Assume that no (pseudo-)\'etale cover of $(X,\emptyset)$ dominates a variety of general type. Is it true that $(X,Z)$ admits no pseudo-\'etale cover dominating a pair of log general type? (In other words, with the terminology of Definition \ref{def:ws}, if $X$ is weakly-special, is $X \setminus Z$ also weakly-special?) 
\end{problem}

Theorem \ref{thm:sym_prods_1} and a simple observation on the complement of the big diagonal in the symmetric product of a variety (see Theorem \ref{thm:sym_prods_2}) give a counterexample to the above Puncturing Problems.

\begin{alphtheorem}[Counterexample to Hassett--Tschinkel's arithmetic puncturing problem, proven in Section~\ref{section:potential_density}]\label{thm:ce_ar}
Let $C$ be a smooth proper geometrically connected curve of genus $g\geq 2$ over a number field $K$, and let $m\geq g$. Define $X:=\Sym^m(C\times \mathbb{P}^1_K)$. Then the following statements hold.
\begin{enumerate}
\item There is a finite field extension $L/K$ such that $X(L)$ is dense, i.e., the normal projective variety $X_{\overline{K}}$ is arithmetically-special over $\overline{K}$.
\item Integral points on the pair $(X, Z)$, where $Z$ is the big diagonal, are not potentially dense and $\mathrm{codim}_X(Z) \geq 2$, i.e., the variety $X_{\overline{K}}\setminus Z_{\overline{K}}$ is not arithmetically-special over $\overline{K}$ and $X\setminus Z\subset X$ is a big dense open.
\item The normal projective variety $X$ has canonical singularities.
\end{enumerate}
\end{alphtheorem}

The counterexample to the Arithmetic Puncturing Problem also provides a counterexample to the Geometric Puncturing Problem.

\begin{alphremark}[Counterexample to Hassett--Tschinkel's geometric puncturing problem, proven in Section~\ref{section:ascending_descending}]\label{thm:ce_geom}
Let $C$ be a smooth proper connected curve of genus $g\geq 2$ over an algebraically closed field $k$ of characteristic zero, and let $m\geq g$. Define $X := \Sym^m(C\times \mathbb{P}^1_k)$. Then the following statements hold.
\begin{enumerate}
\item No finite \'etale cover of $X$ dominates a positive-dimensional variety of general type.
\item The pair $(X,Z)$, where $Z$ is the big diagonal, has a pseudo-\'etale cover which dominates a pair of log-general type.
\item The normal projective variety $X$ has canonical singularities.
\end{enumerate}
\end{alphremark}
 
 Our counterexamples $\Sym^m(C\times \mathbb{P}^1)$ to the above puncturing problems were already mentioned in \cite[p.384]{CCR}.
In fact, our ``smallest'' example $V = \Sym^2(C \times \mathbb{P}^1)$, with $C$ a smooth projective genus two curve, is a special fourfold which becomes non-special after removing a closed subset of codimension two.

Although our example involves a singular projective variety $X$, we note that a desired application of a positive answer to the Arithmetic Puncturing Problem, namely the potential density of rational points on K3 surfaces \cite[Remark~2.14]{HT}, required a positive answer in the singular context (which turns out to be false). Indeed, our construction and argument (in the arithmetic case) are parallel to those in Hassett and Tschinkel's \cite[Remark~2.14]{HT}, except that they look at $\Sym^n(S)$ for a K3 surface $S$, whereas we consider the case $S = C \times \mathbb{P}^1$.

We point out that the varieties $\Sym^m(C\times \mathbb{P}^1)$ can also be used to answer in the negative a question of Kamenova--Lehn on the behaviour of the Kobayashi pseudometric; see Remark \ref{remark:kamenova}.

Despite the fact that Hassett-Tschinkel's conjectures are false for varieties with canonical singularities, it seems reasonable to suspect that they are true for smooth varieties. In the next section we propose corrected conjectures guided by Campana's conjectures.

\subsection{The corrected puncturing conjectures} 
Our starting point is the following ``puncturing'' property for \emph{smooth} special varieties.

\begin{alphtheorem}[Proven in Section~\ref{section:puncturing_specials}] \label{thm:punc_cam}
Let $X$ be a smooth special variety over $k$, and let $U\subset X$ be a dense open whose complement is of codimension at least two. Then $U$ is special. 
\end{alphtheorem}

Note that this is an example of a purity statement. Other examples of such purity statements include, for example, that the fundamental group of $X$ is isomorphic to the fundamental group of $U$ or that the natural restriction map $\mathrm{Br}(X) \to \mathrm{Br}(U)$ of Brauer groups is an isomorphism. Theorem \ref{thm:punc_cam} fails without smoothness assumptions as we have illustrated using $\Sym^m(C\times \mathbb{P}^1)$ (see Remark~\ref{thm:ce_geom}), and so do the purity statements for $\pi_1$ and $\mathrm{Br}$.

Campana's conjectures (Conjecture \ref{conj:campana_analytic} and Conjecture \ref{conj:arithmetic}) combined with Theorem \ref{thm:punc_cam} thus predict that every notion of specialness for a smooth variety is preserved after passing to an open whose complement is of codimension at least two. 

The following result fits in well with the above prediction; it is a consequence of the classical theorem on the invariance of Kobayashi's pseudometric on a smooth variety deprived of a closed subset of codimension at least two \cite[Theorem~3.2.19]{KobayashiBook} (see \cite{CampbellHowardOchiai, CampbellOgawa}).

\begin{theorem}[Campbell--Ogawa, Campbell--Howard--Ochiai]\label{thm:punc_koba}
Let $X$ be a smooth Kobayashi-special variety over $\mathbb{C}$, and let $U\subset X$ be a dense open whose complement is of codimension at least two. Then $U$ is Kobayashi-special.
\end{theorem}

In the case of Brody-specialness, arithmetic-specialness, and geometric-specialness, the expected puncturing property is not known. This leads to the following conjecture. 
 
\begin{conjecture}[The puncturing conjectures]\label{conj:punc}
Let $X$ be a smooth variety over $k$, and let $Z\subset X$ be a closed subset of codimension at least two. Then the following statements hold.
\begin{enumerate}
\item If $k=\CC$ and $X$ is Brody-special, then $X \setminus Z$ is Brody-special. 
\item If $X$ is geometrically-special over $k$, then $X \setminus Z$ is geometrically-special over $k$.
\item If $X$ is arithmetically-special over $k$, then $X \setminus Z$ is arithmetically-special over $k$. 
\item If $X$ has the arithmetic weak Hilbert property over $k$, then $X\setminus Z$ has the arithmetic weak Hilbert property over $k$. 
\end{enumerate}
\end{conjecture}

Note that Conjecture \ref{conj:punc} is similar to the Puncturing Problems of Hassett-Tschinkel, but with four important differences: 
\begin{itemize}
\item we restrict to \emph{smooth} varieties, 
\item we allow $X$ to be \emph{non-proper}, 
\item we propose \emph{additional} conjectures for Brody-special and geometrically-special varieties as well as for varieties satisfying the potential weak Hilbert property,
\item we replace ``weakly-special'' by ``special''. (In this paper, we ignore the question of whether a smooth weakly-special variety remains weakly-special after puncturing.)
\end{itemize} 

Let us discuss some supporting evidence for Conjecture \ref{conj:punc}. For example, as rationally connected varieties are special \cite[Corollary~2.28]{Ca04}, it is natural to study Conjecture \ref{conj:punc} for such varieties. Campana--Winkelmann showed that complements of small closed subsets in a smooth projective rationally connected variety admit a dense entire curve (hence are Brody-special); see \cite{CampanaWinkelmann}. Prior to their work it was not even known whether all rationally connected varieties admit a dense entire curve. On the other hand, since we do not know whether every rationally connected smooth projective variety (or even every smooth projective Fano variety) is arithmetically-special, we also do not know this for complements of small closed subsets in such varieties, except in some special cases \cite{McKinnonRoth, McKinnonZhu}. On the positive side, it is not hard to verify that rationally connected smooth varieties are geometrically-special \cite[Proposition~2.14]{JR}, and that such varieties remain rationally connected (hence geometrically-special) after removing a closed subset of codimension at least two.

Now, for $A$ an abelian variety and $Z$ a closed subset of codimension at least two, since $A$ is special, the variety $A\setminus Z$ is special (Theorem \ref{thm:punc_cam}). It is thus natural to test Conjecture \ref{conj:punc} for abelian varieties. The existence of a dense entire curve in $A$ is a consequence of the fact that it is uniformised by affine space (see \cite[Proposition~3.3]{JR}), i.e., abelian varieties are Brody-special. A proof of the fact that the complement of a small closed subset of an abelian variety is (still) Brody-special was given by Vojta \cite[Proposition~3.2]{VojtaExc}. On the arithmetic side, it is well-known that abelian varieties are arithmetically-special by Frey--Jarden's work on abelian varieties \cite{FreyJarden}. However, proving the arithmetic specialness of $A\setminus Z$ is a notoriously hard problem; it can be verified if $A$ is a product of elliptic curves or if $Z$ consists of the origin and $A$ is a simple CM abelian variety; see \cite[Example~4.4]{HT}. Heuristics motivated by the Arithmetic Puncturing Problem are given in \cite{Siksek} and \cite{KrTsc}. Thus, the arithmetic picture remains essentially completely unresolved (even for abelian surfaces). On the positive side, in the analogous (isotrivial) function field setting, one can prove the geometric specialness of $A \setminus Z$ for any closed subset $Z\subset A$ of codimension at least two in a complex abelian variety (see \cite{Bartsch}).

Finally, if $G$ is a connected linear algebraic group over $k$, then it is not hard to see that $G$ is special. Let $Z\subset G$ be a closed subset of codimension at least two. Recently, it was shown that $G\setminus Z$ is geometrically-special \cite{Bartsch}, and in \cite{Luger4} it was shown by Luger that $G\setminus Z$ satisfies the arithmetic weak Hilbert property (and hence is arithmetically-special). The proof of the arithmetic weak-Hilbert property for $G\setminus Z$ uses strong approximation for semisimple simply connected algebraic groups, and that big opens in such groups still satisfy a form of strong approximation; this form of ``purity'' for smooth varieties with strong approximation was asked about by Wittenberg \cite[Question~2.11]{Wittenberg18}.

\begin{ack} 
The first-named author thanks Jonas Ehrhard for a helpful discussion on Lemma \ref{lemma:gendiagonal}. 
The second-named author gratefully acknowledges J\"org Winkelmann for explaining the proof of Proposition \ref{prop:cw}.(3). We are grateful to Fr\'ed\'eric Campana and Erwan Rousseau for many helpful discussions on special varieties. We are grateful to Daniel Loughran for helpful discussions on the proof of Theorem \ref{thm:hp_for_twists_of_powers_of_P1}.  We thank Olivier Wittenberg for helpful comments and Remark \ref{remark:wittenberg}. 
The second-named author gratefully acknowledges support from the IHES. The third-named author was supported in part by NSF grants DMS-2001205 and DMS-2302298, and a Simons Fellowship from the Simons Foundation.
\end{ack}

\section{Campana's special varieties}\label{section:puncturing_specials}

Let $k$ be an algebraically closed field of characteristic zero. Let $(X,D)$ and $(X',D')$ be snc pairs over $k$. A \emph{morphism $(X,D)\to (X',D')$ of snc pairs} is a morphism $f\colon X\to X'$ such that $f^{-1}(D') \subset D$. A \emph{rational map $(X,D)\ratmap (X',D')$} is a strict rational map $X \setminus D \ratmap X' \setminus D'$, i.e., there is a proper birational surjective morphism $Y \to X \setminus D$ such that the rational map $Y\to X\setminus D \ratmap X'\setminus D'$ is a morphism.

Note that, if $f \colon (X,D) \to (X',D')$ is a morphism of snc pairs, then the morphism $f^\ast \Omega^p_{X'} \to \Omega^p_X$ induces a morphism $f^\ast \Omega^p_{X'}(\log D')\to \Omega^p_X(\log D)$. It suffices to prove this for $p=1$ in which case it is not hard to show \cite[Proposition~11.2]{Iitaka}. 

\begin{lemma}\label{lemma:independence_of_choices}
Let $f \colon (X,D) \to (X',D')$ be a morphism of snc pairs such that $X \setminus D \to X' \setminus D'$ is proper birational. Let $1 \leq p \leq \dim X$ be an integer. Then, $(X,D)$ has a Bogomolov sheaf of rank $p$ if and only if $(X', D')$ has one.
\end{lemma}
\begin{proof}
If $\mathcal{L}'$ is a Bogomolov sheaf of rank $p$ on $(X', D')$, then a nonzero morphism $\mathcal{L'} \to \Omega^p_{X'}(\log D')$ induces a nonzero morphism $f^*\mathcal{L}' \to f^*\Omega^p_{X'}(\log D')$ and via composition with the natural pullback map $f^* \Omega^p_{X'}(\log D') \to \Omega^p_X(\log D)$ we obtain a nonzero morphism $f^*\mathcal{L}' \to \Omega^p_X(\log D)$. As we have $\kappa(\mathcal{L}') = \kappa(f^*\mathcal{L}')$, we see that $f^*\mathcal{L}'$ is a Bogomolov sheaf of rank $p$ on $(X, D)$.

Conversely, let $\mathcal{L}$ be a Bogomolov sheaf of rank $p$ on $(X, D)$. Let $U' \subseteq X'$ denote the maximal open over which $f$ is an isomorphism and let $U := f^{-1}(U')$. Note that the complement of $U'$ in $X'$ has codimension at least two. Then $f$ identifies $D \cap U$ with $D' \cap U'$, so that $(f_* \Omega^p_X(\log D))|_{U'} = \Omega^p_{X'}(\log D')|_{U'}$. Thus, $(f_* \mathcal{L})|_{U'}$ admits a nonzero morphism to $(\Omega^p_{X'}(\log D'))|_{U'}$. As $X'$ is a smooth variety, the line bundle $(f_*\mathcal{L})|_{U'}$ on $U'$ extends to a line bundle $\widetilde{\mathcal{L}}$ on $X'$ and by Hartogs' Lemma, the morphism of sheaves $(f_*\mathcal{L})|_{U'} \to (\Omega^p_{X'}(\log D'))|_{U'}$ extends to a nonzero morphism $\widetilde{\mathcal{L}} \to \Omega^p_{X'}(\log D')$. By construction, we have $\widetilde{\mathcal{L}}(X') = (f_*\mathcal{L})(U') = \mathcal{L}(U)$ (and similarly for tensor powers of $\widetilde{\mathcal{L}}$), so that $\mathcal{L}(X) \subseteq \widetilde{\mathcal{L}}(X')$ and consequently $\kappa(\widetilde{\mathcal{L}}) \geq \kappa(\mathcal{L}) = p$. Hence, $\widetilde{\mathcal{L}}$ is a Bogomolov sheaf of rank $p$ on $(X',D')$, as desired.
\end{proof}

\begin{lemma} \label{lemma:pullback_bs}
Let $f \colon (X,D) \to (X',D')$ be a surjective morphism of snc pairs. If $\mathcal{L}$ is a Bogomolov sheaf for $(X',D')$, then $f^\ast \mathcal{L}$ is a Bogomolov sheaf for $(X,D)$.
\end{lemma}
\begin{proof}
By definition, there is an integer $p$ such that $\mathcal{L}$ admits a nonzero morphism to $\Omega^p_{X'}(\log D')$ and such that $\kappa(\mathcal{L}) = p$. Since $f$ is surjective (and separable), the morphism $f^\ast \Omega^p_{X'}(\log D')\to \Omega^p_X(\log D)$ is injective \cite[Proposition~11.2]{Iitaka}. In particular, the line bundle $f^\ast \mathcal{L}$ admits a nonzero morphism to $\Omega^p_X(\log D)$. Lastly, note that $\kappa(f^*\mathcal{L}) = \kappa(\mathcal{L})$, which finishes the proof. 
\end{proof}

\begin{remark}\label{remark:campana_def} 
Our definition of a special snc pair (Definition \ref{def:special_pair}) coincides with Campana's definition \cite[D\'{e}finition~5.17]{Ca11} by \cite[Th\'eor\`eme~9.9]{Ca11}.
Indeed, while Campana defines a Bogomolov sheaf as being a subsheaf $\mathcal{L} \subseteq \Omega^p_X$, the sections of $\mathcal{L}^{\otimes n}$ he considers to define his $\kappa(X|D, \mathcal{L})$ are really the global sections of the saturations $(\mathcal{L}^{\otimes n})^{\sat} \subseteq \Sym^n(\Omega^p_X(\log D))$, cf.\ the paragraph ``Notations'' at the beginning of \cite[Section~3]{Ca11}. As a result, since $(\mathcal{L}^{\sat})^{\otimes n} = (\mathcal{L}^{\otimes n})^{\sat}$, our Bogomolov sheaves are the saturation in $\Omega^p_X(\log D)$ of the Bogomolov sheaves discussed in \cite{Ca11}.
\end{remark}

A generically finite morphism $(X',D') \to (X,D)$ of snc pairs is an \emph{\'etale covering} if $X'\setminus D'\to X \setminus D$ is finite \'etale.
It is a highly non-trivial fact that a finite \'etale cover of a special snc pair is special.  

\begin{theorem}[Campana]\label{thm:campana_cw}
Let $(X,D)$ be a special snc pair. Let $(X',D')\to (X,D)$ be an \'etale covering. Then $(X',D')$ is special. 
\end{theorem}
\begin{proof}
By Remark \ref{remark:campana_def}, we may appeal to Campana's theorem \cite[Th\'eor\`eme~10.11]{Ca11}.
\end{proof}

An snc pair $(X,D)$ is of \emph{general type} if $\omega_X(D)$ is a big line bundle on $X$. We note Campana's observation that a special snc pair does not dominate a positive-dimensional snc pair of general type (this follows from the far more general \cite[Th\'eor\`eme~9.9]{Ca11}).

\begin{proposition}[Campana]\label{prop:almost_ws}
Let $f \colon (X,D) \to (X',D')$ be a dominant rational map of snc pairs, where $(X,D)$ is special. If $(X',D')$ is of general type, then $\dim X'=0$.   
\end{proposition}

\begin{definition} 
An snc pair $(X,D)$ is \emph{weakly-special} if, for every \'etale covering $(X',D')\to (X,D)$, the snc pair $(X',D')$ does not admit a dominant rational map $(X',D') \to (Z,D_Z)$ to an snc pair of general type $(Z,D_Z)$ with $\dim Z>0$. 
\end{definition}

Combining Theorem \ref{thm:campana_cw} and Proposition \ref{prop:almost_ws} gives the following result of Campana.
 
\begin{corollary}[Campana]\label{corollary:special_is_weakly_special}
If $(X,D)$ is a special snc pair, then $(X,D)$ is weakly-special. 
\end{corollary}


\subsection{Puncturing, images, and birational invariance}
With the Bogomolov sheaf-theoretic definition of a special variety, the fact that smooth special varieties remain special after puncturing is not difficult:
 
\begin{proof}[Proof of Theorem \ref{thm:punc_cam}]
Let $(\overline{X},D)$ be an snc compactification of $X$ and denote by $Z \subseteq \overline{X}$ the closure of $X \setminus U$ in $\overline{X}$. Let $\psi \colon X' \to \overline{X}$ be a proper birational surjective morphism which is an isomorphism over $\overline{X} \setminus Z$ such that $E := X' \setminus \psi^{-1}(U)$ is an snc divisor. (Thus, $(X',E)$ is an snc compactification of $U$.)
To prove the theorem, we have to show that the snc pair $(X',E)$ is special. Indeed, assume that $(X',E)$ were not special. Then, there is an integer $p \geq 1$ and a Bogomolov sheaf $\mathcal{L}' \subseteq \Omega^p_{X'}(\log E)$. Consider the pushforward sheaf $\psi_* \mathcal{L}'$ on $\overline{X}$. As $\psi$ is an isomorphism over the open subset $\overline{X} \setminus Z$, the restriction $(\psi_* \mathcal{L}')|_{\overline{X} \setminus Z}$ is a line bundle on $\overline{X} \setminus Z$. Moreover, we have that $(\psi_* \mathcal{L'})|_{\overline{X} \setminus Z} \subseteq \Omega^p_{\overline{X}}(\log D)|_{\overline{X} \setminus Z}$. As $U \subseteq X$ has a complement of codimension at least two by assumption, the closed subset $Z \subseteq \overline{X}$ is of codimension at least two. Thus, as $\overline{X}$ is smooth, the line bundle $(\psi_* \mathcal{L}')|_{\overline{X} \setminus Z}$ on $\overline{X} \setminus Z$ extends to a sub-line bundle $\mathcal{L} \subseteq \Omega^p_{\overline{X}}(\log D)$ on $\overline{X}$ by Hartogs' Lemma. Now observe that by construction, we have $(\mathcal{L}')^{\otimes n}(X' \setminus E) = (\mathcal{L}^{\otimes n})(\overline{X})$ for every integer $n$. Hence, we have inclusions $(\mathcal{L}')^{\otimes n}(X') \subseteq \mathcal{L}^{\otimes n}(\overline{X})$. This shows that the Iitaka dimensions of these line bundles satisfy $\kappa(\mathcal{L}) \geq \kappa(\mathcal{L}')$. Consequently, $\mathcal{L}$ is a Bogomolov sheaf on $X$, contradicting our assumption that $X$ is special. So we see that $\mathcal{L'}$ cannot exist, so that the pair $(X',E)$ has no Bogomolov sheaves and is hence special.
\end{proof}

\begin{remark}
The assumption in Theorem \ref{thm:punc_cam} on the codimension is obviously necessary. Indeed, $\mathbb{G}_m$ is special, however $\mathbb{G}_m\setminus \{1\}$ is (hyperbolic and) not special.
\end{remark}
 
Lemma \ref{lemma:independence_of_choices}, Lemma \ref{lemma:pullback_bs} and Theorem \ref{thm:campana_cw} imply the following basic properties of special varieties.
 
\begin{lemma} \label{lemma:image_of_special}
Let $f\colon X\to Y$ be a surjective morphism of varieties. Then the following statements hold.
\begin{enumerate}
\item If $X$ is special, then $Y$ is special.
\item If $f$ is proper and birational, then $X$ is special if and only if $Y$ is special.
\item If $f$ is finite \'etale, then $X$ is special if and only if $Y$ is special.
\end{enumerate}
\end{lemma}

A notion closely related to specialness is that of a weakly-special variety:

\begin{definition}\label{def:ws}
We say that a variety $X$ is \emph{weakly-special} if there is a resolution of singularities $X'\to X$ and an snc compactification $\overline{X'}$ of $X'$ with boundary $D$ such that $(\overline{X'},D)$ is weakly-special. 
\end{definition}

If $X$ is a variety, then $X$ is weakly-special if and only if, for every resolution of singularities $X'\to X$ and every snc compactification $\overline{X'}$ of $X'$ with $D := \overline{X'}\setminus X'$, the pair $(\overline{X'},D)$ is weakly-special. In other words, the notion of being weakly-special is independent of the choice of snc model. 
  
We note that our definition of a weakly special variety differs from that of Cadorel--Deng--Yamanoi \cite[Definition~1.10]{CadorelDengYamanoi} for singular varieties.
Essentially, this is because a singular model of a variety may have fewer \'etale covers than its smooth models.  
  More precisely, note that Cadorel--Deng--Yamanoi begin with a normal variety $X$ (possibly singular) and require that no finite \'etale cover of $X$ rationally dominates a positive-dimensional variety of general type. 
In contrast, we impose this condition  after passing to a resolution $\widetilde{X}$ of $X$. Since the induced homomorphism on fundamental groups $\pi_1(\widetilde{X}) \to \pi_1(X)$ is   not necessarily an isomorphism, the two notions can differ. 

For an explicit example, let $C$ be a smooth projective curve of genus two, and let $E$ be an elliptic curve.
Let $\sigma$ be the hyperelliptic involution on $C$ and let $\tau$ be the translation by a point of order two on $E$.
Then $(\sigma, \tau)$ defines an involution on $C \times E$ and we let $Y = (C \times E)/\langle (\sigma, \tau)\rangle$ be the quotient.
Although $Y$ does not dominate any positive-dimensional variety of general type, it is not weakly special, since $C \times E \to Y$ is finite \'etale and $C \times E$ dominates $C$.
Choose an embedding $Y \subseteq \mathbb{P}^n$, and let $X \subseteq \mathbb{P}^{n+1}$ denote the projective cone over $Y$.
Let $X' \to X$ be any resolution of singularities.
We claim that $X$ is not weakly special in the sense used above, yet it is weakly special in the sense of \cite[Definition~1.10]{CadorelDengYamanoi}.
Indeed, this can be seen from the following observations.
\begin{enumerate}
\item The threefold $X'$ admits a finite \'etale cover that dominates $C$, and is therefore not weakly special. Consequently, the threefold $X$ is also not weakly special in the sense defined above.
\item The threefold $X$ admits no non-trivial finite \'etale covers (see \cite[Exercise~4.8.1]{DebarreBook}).
\item The threefold $X'$ does not admit a dominant rational map to a positive-dimensional variety $Z$ of general type. Indeed, such a map $X' \ratmap Z$ would factor over $Y$, contradicting the fact that $Y$ does not dominate any positive-dimensional variety of general type.  
We thus conclude that $X$ is weakly special in the sense of \cite[Definition~1.10]{CadorelDengYamanoi}.
\end{enumerate} 

The following corollary due to Campana follows directly from the definitions and Corollary \ref{corollary:special_is_weakly_special}. 
 
\begin{corollary}[Campana]\label{cor:sp_is_ws}
If $X$ is a special variety, then $X$ is weakly-special. 
\end{corollary}

\begin{remark}  
If $X$ is proper and $\dim X \leq 2$, then the converse to Corollary \ref{cor:sp_is_ws} holds. Indeed, this is trivial for curves and for surfaces follows by going through the Enriques--Kodaira classification (see \cite[Corollary~3.33]{Ca04} for a classification of special surfaces). If $\dim X \geq 3$, there are examples of weakly-special non-special smooth projective varieties; see \cite{BT, RTW, BCJW}.
\end{remark}

\section{Ascending and descending specialness properties} \label{section:ascending_descending}
Before we prove Theorem \ref{thm:sym_prods_2}, we state and prove two well-known lemmas on the class of special varieties.

\begin{proposition}\label{prop:images}
Let $X \to Y$ be a dominant morphism of varieties over $k$. Then the following statements hold.
\begin{enumerate}
\item If $X$ is special, then $Y$ is special.
\item If $X$ is weakly-special, then $Y$ is weakly-special.
\item If $k=\CC$ and $X$ is Brody-special, then $Y$ is Brody-special.
\item If $k=\CC$ and $X$ is Kobayashi-special, then $Y$ is Kobayashi-special.
\item If $X$ is arithmetically-special, then $Y$ is arithmetically-special.
\item If $X$ is geometrically-special, then $Y$ is geometrically-special.
\end{enumerate}
\end{proposition}
\begin{proof}  
If $X$ is special, then $Y$ is special by Lemma \ref{lemma:image_of_special}; this proves $(1)$.
If $Y$ is not weakly-special, then it follows easily from the definition that $X$ is not weakly-special. This proves $(2)$. 
To prove $(3)$, compose a dense entire curve in $X^{\an}$ with the dominant map $X^{\an}\to Y^{\an}$ to obtain a dense entire curve in $Y^{\an}$.
To prove $(4)$, use the distance-decreasing property of the Kobayashi pseudometric. 
To prove $(5)$, use that the image of a dense subset of near-integral points on $X$ along $X\to Y$ is a dense subset of near-integral points on $Y$ (after choosing suitable models over a suitable $\mathbb{Z}$-finitely generated subring of $k$).
Finally, $(6)$ is proven in \cite[Section~2.2]{JR}.
\end{proof}

\begin{proposition}\label{prop:cw}
Let $X \to Y$ be a finite \'etale morphism of (integral) varieties over $k$. Then the following statements hold.
\begin{enumerate}
\item The variety $X$ is special if and only if $Y$ is special.
\item The variety $X$ is weakly-special if and only if $Y$ is weakly-special.
\item If $k=\CC$, then $X$ is Brody-special if and only if $Y$ is Brody-special.
\item If $k=\CC$, then $X$ is Kobayashi-special if and only if $Y$ is Kobayashi-special.
\item The variety $X$ is arithmetically-special if and only if $Y$ is arithmetically-special.
\item The variety $X$ is geometrically-special if and only if $Y$ is geometrically-special.
\end{enumerate}
\end{proposition}
\begin{proof}   
First note that $X \to Y$ is surjective. Thus, if $X$ is special (resp. weakly-special, Brody-special, etc.), then it follows from Proposition \ref{prop:images} that $Y$ is so as well. We now prove the converse statements. 

If $Y$ is special, then $X$ is special by Lemma \ref{lemma:image_of_special}; this proves $(1)$.
Furthermore, it follows directly from the definition that if $Y$ is weakly-special, then $X$ is weakly-special; this proves $(2)$. 
Also, note that $(3)$ follows from the fact that entire curves lift along finite \'etale morphisms. It remains to prove $(4)$, $(5)$ and $(6)$.

To prove $(4)$, assume $Y$ is Kobayashi-special. Let $Y'\to Y$ be a resolution of singularities and let $X'\to X$ be the pullback. Note that $X'$ is a smooth (integral) variety. Since $Y$ is Kobayashi-special, we have that $d_{Y'}\equiv 0$. To show that $X$ is Kobayashi-special, it suffices to show that $d_{X'}\equiv 0$. Therefore, replacing $Y$ by $Y'$ and $X$ by $X'$, we may assume that $X$ and $Y$ are smooth. To show that $d_X\equiv 0$, assume that there are distinct points $p_1$ and $p_2$ in $X$ such that $d_X(p_1,p_2)>0$. Define the equivalence class of a point $P$ in a complex-analytic space $\mathcal{X}$ to be the locus of points $Q$ such that $d_{\mathcal{X}}(P,Q) =0$. Then, since the Kobayashi pseudometric defines a continuous function on $X\times X$, the equivalence class $X_1$ of $p_1$ is closed in $X$. Moreover, since $d_X(p_1,p_2) >0$, this equivalence class is disjoint from the (closed) equivalence class $X_2$ of $p_2$. Moreover, the formula for the pseudometric $d_Y$ given in \cite[Theorem~3.2.8.(1)]{KobayashiBook} shows that the equivalence class of \emph{any} point $P$ of $X$ surjects onto $Y$. Since $X^{\an}$ is connected, we have that $X^{\an} \neq X_1 \cup X_2$. Therefore, there is a point $p_3$ in $X$ such that $d_X(p_1,p_3) >0$ and $d_X(p_2,p_3) >0$. Thus, the equivalence class $X_3$ of $p_3$ is a closed subset disjoint from $X_2$ and $X_3$. If $n:=\deg(X\to Y)$, then repeating this process gives a sequence of closed subsets $X_1,\ldots, X_n$ which are pairwise disjoint. Since the covering $X\to Y$ is of degree $n$, we see that $X= X_1\sqcup \ldots \sqcup X_n$ contradicting the connectivity of $X$. This proves $(4)$. 

Note that $(5)$ is a consequence of a (fairly general) version of the Chevalley-Weil theorem. Due to lack of reference in the near-integral setting we include a proof. We closely follow \cite[Lemma~8.2]{JLitt}. 
Assume that $Y$ is arithmetically-special over $k$. Choose a regular $\mathbb{Z}$-finitely generated integral domain $A\subset k$, a finite type separated model $\mathcal{X}$ for $X$ over $A$, a finite type separated model $\mathcal{Y}$ for $Y$ over $A$, and a finite \'etale surjective morphism $F \colon \mathcal{X}\to \mathcal{Y}$ extending $X\to Y$ such that $\mathcal{Y}(A)^{(1)}$ is dense in $Y$. For every near-integral point $y\in \mathcal{Y}(A)^{(1)}$, there exist a dense open subscheme $U_y\subset \Spec A $ whose complement in $\Spec A$ is of codimension at least two and a morphism $U_y\to \mathcal{Y}$. Pulling back $U_y\to \mathcal{Y}$ along $F \colon \mathcal{X}\to \mathcal{Y}$, we obtain a finite \'etale surjective morphism $V_y:= U_y\times_{\mathcal{Y}} \mathcal{X}\to U_y$ of degree $\deg(f)$ which, by purity of the branch locus extends to a finite \'etale morphism $\overline{V_y}\to \Spec A$. By the Hermite-Minkowski theorem for arithmetic schemes \cite{smallness}, the set of isomorphism classes of the $V_y$ is finite as $y$ runs over $\mathcal{Y}(A)^{(1)}$. In particular, there is a $\ZZ$-finitely generated integral domain $B\subset k$ containing $A$ such that some dense subset of $\mathcal{Y}(A)^{(1)}$ lies in the image of $\mathcal{X}(B)^{(1)}$. This implies that the latter is dense, as required.
 
Finally, to conclude the proof, note that $(6)$ is \cite[Lemma~2.11]{JR}.
\end{proof}

As an application of the above propositions, we make the simple observation that if the complement of the big diagonal in a symmetric power of $X$ is special, then $X$ is forced to be special. We also prove the analogous statement for every other notion of specialness. 
 
\begin{theorem}\label{thm:sym_prods_2}
Let $X$ be a variety over $k$, let $n\geq 1$ be an integer, and let $Z\subset \Sym^n(X)$ be the big diagonal. Then the following statements hold.
\begin{enumerate}
\item If $X$ is not special, then $\Sym^n(X)\setminus Z$ is not special.
\item If $X$ is not weakly-special, then $\Sym^n(X)\setminus Z$ is not weakly-special.
\item If $k=\CC$ and $X$ is not Brody-special, then $\Sym^n(X)\setminus Z$ is not Brody-special.
\item If $k=\CC$ and $X$ is not Kobayashi-special, then $\Sym^n(X)\setminus Z$ is not Kobayashi-special.
\item If $X$ is not arithmetically-special over $k$, then $\Sym^n(X)\setminus Z$ is not arithmetically-special over $k$.
\item If $X$ is not geometrically-special over $k$, then $\Sym^n(X)\setminus Z$ is not geometrically-special over $k$.
\end{enumerate}
\end{theorem}
\begin{proof} 
Note that $X^n\setminus \Delta\to \Sym^n(X)\setminus Z$ is finite \'etale. Thus, if $\Sym^n(X)\setminus Z$ is special, then $X^n \setminus \Delta$ is special (Proposition \ref{prop:cw}). Now, since the special variety $X^n\setminus \Delta$ surjects onto $X$ (use the composition of the inclusion $X^n \setminus \Delta \subset X^n$ with a projection map $X^n\to X$), it follows from Proposition \ref{prop:images} that $X$ is special. This proves $(1)$. 
 
The same line of reasoning also proves $(2)$, $(3)$, $(4)$, $(5)$ and $(6)$. 
\end{proof}

We can now show that $\Sym^m(C \times \mathbb{P}^1)$ for $C$ a smooth projective curve of genus $g \geq 2$ and $m \geq g$ gives a counterexample to Hassett--Tschinkel's geometric puncturing problem (Problem~\ref{pGP}), that is, we can now prove Remark~\ref{thm:ce_geom}.

\begin{proof}[Proof of Remark~\ref{thm:ce_geom}]
That $\Sym^m(C \times \mathbb{P}^1)$ is weakly-special follows from Theorem~\ref{thm:sym_prods_ccr} and Corollary~\ref{cor:sp_is_ws}; this shows $(1)$.
The complement of the big diagonal in $\Sym^m(C \times \mathbb{P}^1)$ is not weakly-special by Theorem~\ref{thm:sym_prods_2}; this shows $(2)$.
Thus, it remains to show that $\Sym^m(C \times \mathbb{P}^1)$ has canonical singularities. This follows from the fact that the Hilbert scheme $\mathrm{Hilb}^m(C\times \mathbb{P}^1_k)$ of closed subschemes of length $m$ on $C\times \mathbb{P}^1$ provides a crepant resolution of singularities of $\Sym^m(C \times \mathbb{P}^1)$ \cite[Theorem~7.4.6]{BrionKumar}.   
\end{proof}

We finish with a discussion of a question of Kamenova--Lehn \cite[Question~3.7]{KamenovaLehn}.

\begin{remark}\label{remark:kamenova}
Let $C$ be a smooth projective curve of genus $g$ over $\mathbb{C}$. Let $m\geq g$. Then the variety $\Sym^m(C\times \mathbb{P}^1)$ can be used to give a negative answer to a question of Kamenova and Lehn \cite[Question~3.7.(1)]{KamenovaLehn}. Indeed, we know that $\Sym^m(C\times \mathbb{P}^1)$ is Kobayashi-special (Theorem \ref{thm:sym_prods_ccr}) with canonical (hence log-terminal) singularities (Remark~\ref{thm:ce_geom}.(3)). However, the complement of the big diagonal in $\Sym^m(C\times \mathbb{P}^1)$ is smooth and not Kobayashi-special (Theorem \ref{thm:sym_prods_2}.(4)). 
\end{remark}

\section{A criterion for density of graphs} \label{section:density_criteria}

When checking whether a given variety $X$ is geometrically-special, one has to check that the graphs of the morphisms $\phi_i \colon C \to X$ one has written down are actually dense in $C \times X$. As this can be sometimes rather difficult, this subsection is dedicated to establishing a criterion that can be slightly easier to check in practice. 

We start by proving some technical lemmas. For some intuition about the first lemma, consider the case where the polynomial $p=1$ is constant. Then the set $\Delta_{p,n}$ described in the lemma is just the usual diagonal. It is clearly closed as $X$ is a separated scheme. Recall that a \emph{numerical polynomial} is a polynomial $p \in \mathbb{Q}[t]$ such that $p(n) \in \mathbb{Z}$ for every $n \in \mathbb{Z}$.

\begin{lemma}\label{lemma:gendiagonal}
Let $X$ be a projective variety with a fixed ample line bundle $\mathcal{L}$. Let $p \in \mathbb{Q}[t]$ be a numerical polynomial and let $n \in \mathbb{N}$. Then the following subset of $X^n = X \times X \times ... \times X$ is closed:
$$ \Delta_{p,n} = \left\{ (x_1,...,x_n)\in X^n~\bigg|~\begin{aligned} &\{x_1,...,x_n\}~\text{is contained in a closed subscheme of}~X~\\ &\text{with Hilbert polynomial}~p \end{aligned} \right\} $$~
\end{lemma}
\begin{proof}
Consider the Hilbert scheme $H = \Hilb(X,p)$ which parametrizes closed subschemes of $X$ with Hilbert polynomial $p$. Note that $H$ is a projective scheme which comes equipped with a universal family $\mathcal{F} \subseteq X \times H$, which is a closed subscheme of $X \times H$ (and, set-theoretically consists of those points $(x,h) \in X \times H$ satisfying $x \in h$).

Now let $\mathcal{G} \subseteq X^n \times H^n$ be the intersection of $\mathcal{F}^n$ and $X^n \times \Delta$ where $\Delta \subseteq H^n$ denotes the diagonal (which is closed as $H$ is projective). Then $\mathcal{G}$ is the intersection of two closed subschemes and is hence a closed subscheme of $X^n \times H^n$. As $H$ is proper, the projection $X^n \times H^n \to X^n$ is closed. Hence the image of $\mathcal{G}$ in $X^n$ is closed. Now, note that this image is precisely the subset $\Delta_{p,n}$.
\end{proof}

Suppose we are given a collection of points in projective space $\mathbb{P}^n$ and want to figure out whether all of them are contained in some line. Then we can check this by looking at all three-element subsets of the collection. In particular, we can check it without ever looking at infinitely many of them at once. We generalize this idea.

\begin{lemma}\label{lemma:testonfinitesubset}
Let $X$ be a projective variety with a fixed ample line bundle $\mathcal{L}$. Let $p \in \mathbb{Q}[t]$ be a numerical polynomial. Let $(x_i)_{i \in I}$ be a collection of closed points of $X$. Suppose that there is no closed subscheme of $X$ with Hilbert polynomial $p$ containing all the $x_i$. Then there is a finite subset $J \subseteq I$ such that the collection $(x_j)_{j \in J}$ also has this property.
\end{lemma}
\begin{proof}
Consider the Hilbert scheme $H = \Hilb(X,p)$ together with the universal family $\mathcal{F} \subseteq X \times H$. For $i \in I$, let $Z_i \subseteq H$ be the fiber of the projection $\mathcal{F} \to X$ over the point $x_i$. Set-theoretically, $Z_i$ is the closed subset consisting of all points $h \in H$ satisfying $x_i \in h$. The assumption that no closed subscheme of $X$ with Hilbert polynomial $p$ contains all the $x_i$ means that $\bigcap_{i \in I} Z_i$ is empty. Because $H$ is of finite type over a field, it is quasi-compact. This implies that there is a finite subset $J \subseteq I$ such that $\bigcap_{j \in J} Z_j$ is empty. The finite collection $(x_j)_{j \in J}$ now has the desired properties.
\end{proof}

We can now use the lemmas we just proved to study the graphs of morphisms.

\begin{lemma}\label{lemma:smallimageclosed}
Let $X$ be a projective variety with a fixed ample line bundle $\mathcal{L}$. Let $Y$ be a variety and let $(\phi_i \colon Y \to X)_{i \in I}$ be a family of morphisms. Let $p \in \mathbb{Q}[t]$ be a numerical polynomial. Consider the following subset of $Y$:
$$ Z = \left\{ y \in Y~\bigg|~\begin{aligned} &\{ \phi_i(y)~|~i \in I\}\subseteq X~\text{is contained in a} \\ &\text{closed subscheme of}~X~\text{with Hilbert polynomial}~p \end{aligned} \right\} $$
Then $Z$ is closed in $Y$.
\end{lemma}
\begin{proof}
For a finite subset $J \subseteq I$ we define the following morphism:
$$ \Phi_J \colon Y \to X^{\abs{J}} \quad y \mapsto (\phi_j(y))_{j \in J} $$
By Lemma \ref{lemma:testonfinitesubset}, the failure of infinitely many points to lie on a closed subscheme of some fixed Hilbert polynomial can be detected on a finite subset of them. Consequently, we have:
$$ Z = \bigcap_{J \subseteq I~\text{finite}} \Phi_J^{-1}(\Delta_{p,\abs{J}}) $$
where we used the notation of Lemma \ref{lemma:gendiagonal}. By using that lemma, we see that this equality expresses $Z$ as an intersection of closed subsets, so $Z$ is closed.
\end{proof}

\begin{remark}
Combining the previous lemma with the observation that $\mathbb{Q}[t]$ is countable leads to the following corollary: 
Let $X$ be a projective variety and let $Y$ be a variety. Let $(\phi_i \colon Y \to X)_{i \in I}$ be a family of morphisms. Then the following set is a countable union of closed subvarieties of $Y$:
\[ \{ y \in Y~|~\{ \phi_i(y)~|~i \in I\}\subseteq X~\text{is not dense}\} \]
In particular, when working over an uncountable base field and $\dim(Y) > 0$, the complement of this set is either empty or contains uncountably many points.
We will however not use this statement in the sequel as the conclusion is vacuous when working over countable fields.
\end{remark}

We can now prove our desired criterion for testing the density of the graphs of a family of morphisms. Note that if we assume $k$ to be uncountable, the next theorem immediately follows from the previous remark.

\begin{theorem}\label{thm:geomspeccriterion} 
Let $Y$ be a variety and let $X$ be a quasi-projective variety. Let $(\phi_i \colon Y \to X)_{i \in I}$ be a family of morphisms. Suppose that there is a point $y_0 \in Y$ such that $\{ \phi_i(y_0)~|~i \in I\}$ is dense in $X$. Then $S = \bigcup \Gamma_{\phi_i}$ is dense in $Y \times X$.
\end{theorem}
\begin{proof}
We may assume that $X$ is projective. (Indeed, let $\overline{X}$ be a projective compactification of $X$. Then $S$ is dense in $Y \times X$ if and only if it is dense in $Y \times \overline{X}$.) We now, for the rest of the proof, fix a closed immersion of $X$ into projective space. Doing this allows us to talk about Hilbert polynomials of closed subschemes of $X$.

For the sake of contradiction, suppose that $S$ was not dense in $Y \times X$. Then, there is a proper closed subscheme $Z \subsetneq Y \times X$ containing $S$. By generic flatness, the (surjective) projection morphism $Z \to Y$ is flat over a dense open $Y^o \subseteq Y$. Let $Z^o$ denote the preimage of $Y^o$ in $Z$. It is an open subset of $Z$. Since the Hilbert polynomial of the fibers is independent of the fiber for a flat morphism \cite[Theorem III.9.9]{Har}, every fiber of the projection $Z^o \to Y^o$ has the same Hilbert polynomial $p$. Since $Y \times X$ is irreducible, we must have $\dim(Z) < \dim(Y \times X)$. This implies that $p$ has degree smaller than $\dim X$. As $Z$ contains $S$, this means that for every $y \in Y^o$, the set $\{ \phi_i(y)~|~i \in I\}$ (which is the ``fiber'' of $S$ over $y$) is contained in a closed subscheme of $X$ with Hilbert polynomial $p$ (namely the fiber of $Z^o \to Y^o$ over $y$). Consequently, the dense open $Y^o$ is contained in the subset 
\[ \left\{ y \in Y~\bigg|~\begin{aligned} &\{ \phi_i(y)~|~i \in I\}\subseteq X~\text{is contained in a} \\ &\text{closed subscheme of}~X~\text{with Hilbert polynomial}~p \end{aligned} \right\}. \]
However, by Lemma \ref{lemma:smallimageclosed}, the latter subset is closed. Since it does not contain $y_0$ by assumption, this is a contradiction. So $Z$ cannot exist and we are done.
\end{proof}

\begin{remark}
We can also rephrase Theorem  \ref{thm:geomspeccriterion} as follows: Suppose that $(\phi_i \colon Y \to X)_{i \in I}$ is a family of morphisms from the variety $Y$ to the quasi-projective variety $X$. Consider the induced morphism $Y \times I \to Y \times X$ which sends $(y,i)$ to $(y,\phi_i(y))$ (where we consider the set $I$ as the $I$-indexed disjoint union of copies of $\Spec(k)$). Then, if the restriction $\{y\} \times I \to \{y\} \times X$ is dominant for one point $y \in Y$, the morphism $Y \times I \to Y \times X$ is dominant as well.
\end{remark}

\begin{corollary}\label{cor:geomspeccriterion_arithmetic}
Let $K$ be a field of characteristic zero. Let $Y$ be a variety over $K$ with $Y(K)$ dense and let $X$ be a quasi-projective variety over $K$. Let $(\phi_i \colon Y \to X)_{i \in I}$ be a family of morphisms over $K$. Suppose that there is a point $y_0 \in Y(K)$ such that $\{ \phi_i(y_0)~|~i \in I\}$ is dense in $X$. Then $S = \bigcup \Gamma_{\phi_i}(K)$ is dense in $Y \times X$.
\end{corollary}
\begin{proof}
Since $Y(K)$ is dense in $Y$, for every $i$, we have that $\Gamma_{\phi_i}(K)$ is dense in $\Gamma_{\phi_i}$. Thus, the subset $\bigcup_{i \in I} \Gamma_{\phi_i}(K)$ is dense in $\bigcup_{i \in I} \Gamma_{\phi_i}$. However, the latter is dense in $Y\times X$ by Theorem \ref{thm:geomspeccriterion}.
\end{proof}

\begin{remark} 
In this paper we will use Theorem \ref{thm:geomspeccriterion} to prove that certain symmetric powers are geometrically-special (see Theorems \ref{thm:sym_n_is_geom_special_0} and \ref{thm:symCcrossE_geom_special} below). We will also use Theorem \ref{thm:geomspeccriterion} (or rather its consequence Corollary \ref{cor:geomspeccriterion_arithmetic}) to prove that certain symmetric powers are arithmetically-special (see Theorem \ref{thm:symCcrossE_arith_special}). Finally, in \cite{Bartsch}, Theorem \ref{thm:geomspeccriterion} is used to prove the geometric specialness of every algebraic group.
\end{remark}

\section{Geometrically-special varieties: density of pointed curves}

We recall the definition of a geometrically-special variety (\cite[Definition~1.7]{JR} or Definition \ref{def:geom_specialness}). Throughout this section, $k$ denotes an algebraically closed field of characteristic zero.

\begin{definition}
Let $X$ be a variety over $k$. We say that $X$ is \emph{geometrically-special (over $k$)} if there is a dense subset $S \subseteq X(k)$ such that for every $s \in S$ there is a smooth quasi-projective curve $C$, a closed point $c \in C$ and a family of morphisms $(\phi_i \colon C \to X)_{i \in I}$ satisfying $\phi_i(c) = s$ such that $\bigcup_{i \in I} \Gamma_{\phi_i} \subseteq C \times X$ is Zariski-dense. Here, $\Gamma_{\phi_i}$ denotes the graph of the morphism $\phi_i$.
\end{definition}

Slightly abusing the language, we will call a family of morphisms $(\phi_i)_{i \in I}$ as in the above definition a \emph{covering set for $X$ through $s$}, even though the graphs really cover the product space $C \times X$.
If $X$ is a variety and $U \subseteq X$ is an open subvariety such that $U$ is geometrically-special, then the variety $X$ is geometrically-special as well, as we can simply postcompose any given covering set with the inclusion map $U \to X$.

\subsection{Symmetric powers} 
In this subsection we show that for a curve $C$ of genus $g$, the symmetric powers $\Sym^m(C)$ and $\Sym^m(C \times \mathbb{P}^1)$ are geometrically-special for all $m \geq g$, thereby proving part (5) of Theorem \ref{thm:sym_prods_1}.
We first note the following general lemma.
\begin{lemma}\label{lemma:projbundlegeomspec}  
Let $X$ be a proper variety and let $\mathcal{F}$ be a coherent sheaf on $X$ such that $\mathbb{P}(\mathcal{F})$ is integral. Then there is a dense open $U \subseteq X$ such that, for every $x \in U(k)$, every covering set $(\phi_i \colon (C,c) \to (X,x))_{i\in I}$ and every $y \in \mathbb{P}(\mathcal{F})$ lying over $x$, there is a covering set $(\psi_j \colon (C,c) \to (\mathbb{P}(\mathcal{F}),y))_{j\in J}$. 
\end{lemma}
\begin{proof}
Let $U \subseteq X$ be a nonempty open subscheme over which $\mathcal{F}$ is free. Then there is a natural number $n \geq 0$ such that $\mathcal{F}|_U \cong \mathcal{O}_U^n$. Consequently, we have $\mathbb{P}(\mathcal{F}|_U) \cong \mathbb{P}_U^{n-1}$ as schemes over $U$. This implies in particular that the proper varieties $\mathbb{P}(\mathcal{F})$ and $\mathbb{P}^{n-1} \times X$ are birational. Let $y \in \mathbb{P}(\mathcal{F}|_U)$ lying over a point $x \in U(k)$ for which there is a covering set $(\phi_i \colon (C,c)\to (X,x))_{i \in I}$. We may view $y$ as a point on $\mathbb{P}^{n-1}\times U\subset \mathbb{P}^{n-1}\times X$ and write $y =(y_1,x)$. Let $f \colon (C,c)\to (\mathbb{P}^{n-1},y)$ be a non-constant morphism. Note that the automorphism group $G$ of $(\mathbb{P}^{n-1},y_1)$ acts transitively on $\mathbb{P}^{n-1}\setminus \{y_1\}$. Thus, the collection of morphisms $((g \circ f, \phi_i) \colon (C,c) \to (\mathbb{P}^{n-1}\times X, y))_{g \in G, i \in I}$ is a covering set. Now let $\sigma \colon \mathbb{P}^{n-1} \times X \ratmap \mathbb{P}(\mathcal{F})$ be the birational map induced by the identification of $\mathbb{P}^{n-1} \times U$ and $\mathbb{P}(\mathcal{F})$. Then, by construction, $y$ lies in the regular locus of $\sigma$. Thus, we obtain, for every $i \in I, g \in G$, a rational map $\sigma \circ (g \circ f, \phi_i) \colon (C,c) \ratmap (\mathbb{P}(\mathcal{F}),y)$. Since $C$ is a smooth curve and $\mathbb{P}(\mathcal{F})$ is a proper variety, these rational maps define morphisms. Thus, we obtain a covering set $(\sigma \circ (g \circ f, \phi_i) \colon (C,c) \to (\mathbb{P}(\mathcal{F}),y)_{g \in G, i \in I}$, as desired. 
\end{proof}

The relevance of the lemma for our purposes comes from the following well-known fact; see \cite[Theorem~4]{Schwarzenberger}.

\begin{lemma}\label{lemma:sym_n_coherent_sheaf} 
If $C$ is a smooth projective curve of genus $g_C$ and $m \geq 1$ is an integer, then there is a coherent sheaf $\mathcal{E}$ on $\Pic^m(C)$ and an isomorphism $\Sym^m(C) \cong \mathbb{P}(\mathcal{E})$ of schemes over $\Pic^m(C)$, where the morphism $\Sym^m(C) \to \Pic^m(C)$ comes from viewing an element $[(c_1,...,c_m)] \in \Sym^m(C)$ as the divisor $c_1+...+c_m$ on $C$. If $m \geq g_C$, then the support of $\mathcal{E}$ equals $\Pic^m(C)$, so that $\Sym^m(C)$ is birational to $\mathbb{P}^{m-g_C}\times \Pic^m(C)$. 
\end{lemma}

In the setting of the previous lemma, we see that for dimension reasons, the morphism $\Sym^m(C) \to \Pic^m(C)$ can only be surjective if $m \geq g_C$.
In particular, if $m < g_C$, then the image of $\Sym^m(C) \to \Pic^m(C)$ is a proper closed subvariety of the abelian variety $\Pic^m(C)$.
For $m \geq 1$, this image is not contained in any proper abelian subvariety of $\Pic^m(C)$.
Thus, \cite[Theorem~3.5]{JR} shows that if $1 \leq m < g_C$, the image of $\Sym^m(C)$ in $\Pic^m(C)$ is not geometrically-special and so, $m \geq g_C$ is a necessary condition for $\Sym^m(C)$ being geometrically-special.
As we show now, this condition is also sufficient.

\begin{corollary}\label{cor:sympowergeomspec}
Let $C$ be a smooth projective curve and let $m \geq 1$ be an integer. If $m \geq g_C$, there is a dense open subset $U \subseteq \Sym^m(C)$, a smooth projective curve $D$, and a point $d \in D(k)$, such that for every $u \in U(k)$, there is a covering set $(\psi_j \colon (D,d) \to (\Sym^m(C),u))_{j \in J}$. In particular, the variety $\Sym^m(C)$ is geometrically-special if and only if $m \geq g_C$.
\end{corollary}
\begin{proof}
Note that $\Pic^m(C)$ is an abelian variety, isomorphic to the Jacobian of $C$. In particular, it is projective and geometrically-special \cite[Proposition~3.1]{JR}. Let $o \in \Pic^m(C)$ be any point through which there is a covering set $(\phi_i \colon (D,d) \to (\Pic^m(C),o))_{i \in I}$. Then, as the automorphism group of an abelian variety acts transitively, we see that for any given point $x \in \Pic^m(C)$, there is a covering set $((D,d) \to (\Pic^m(C),x))$. By Lemma \ref{lemma:sym_n_coherent_sheaf}, there is a coherent sheaf $\mathcal{E}$ on $\Pic^m(C)$ such that $\Sym^m(C) \cong \mathbb{P}(\mathcal{E})$. Thus, by Lemma \ref{lemma:projbundlegeomspec}, there is a dense open subset $V \subseteq \Pic^m(C)$ such that there is a covering set from the pointed curve $(D,d)$ through every point of $\Sym^m(C)$ lying over $V \subseteq \Pic^m(C)$. Now let $U$ be the preimage of $V$ in $\Sym^m(C)$ and note that $U$ is nonempty (hence a dense open) as the map $\Sym^m(C) \to \Pic^m(C)$ is surjective. This concludes the proof.
\end{proof}

We will need the following Lemmas in our proof that $\Sym^m(C \times \mathbb{P}^1)$ is geometrically-special.

\begin{lemma}\label{lemma:finitereflectsdensity}
Let $X$, $Y$ be two varieties and let $\phi \colon X \to Y$ be a finite morphism. Let $S \subseteq X(k)$ be a subset. Then, if $\phi(S) \subseteq Y$ is dense, so is $S \subseteq X$.
\end{lemma}
\begin{proof}
Finite morphisms are closed. For closed continuous maps between topological spaces, we have $\phi(\overline{S}) = \overline{\phi(S)}$. Thus $\phi(\overline{S}) = Y$ and in particular $\phi$ is surjective. As finite surjective morphisms preserve dimension, it follows that $\dim(X) = \dim(Y) = \dim(\overline{S})$. As $X$ is irreducible this means $\overline{S} = X$ and we are done.
\end{proof}

\begin{lemma}\label{lemma:finitereflectsdensity2}
Let $X$, $Y$ be two varieties and let $D$ be a curve. Let $\phi \colon X \to Y$ be a finite morphism and let $(\psi_i \colon D \to X)_{i \in I}$ be a family of morphisms. Then, if $\bigcup_{i \in I} \Gamma_{\phi \circ \psi_i} \subseteq D \times Y$ is dense, so is $\bigcup_{i \in I} \Gamma_{\psi_i} \subseteq D \times X$.
\end{lemma}
\begin{proof}
Note that, for every $i \in I$, we have $(\id_D,\phi)(\Gamma_{\psi_i}) = \Gamma_{\phi \circ \psi_i}$. Thus we conclude by Lemma \ref{lemma:finitereflectsdensity}.
\end{proof}

\begin{lemma}\label{lemma:interpolationpolys}
Let $(x_1,...,x_m)$ be an $m$-tuple of pairwise distinct closed points of $\mathbb{P}^1$. Let $(y_1,...,y_m)$ be any other $m$-tuple of closed points of $\mathbb{P}^1$. Then there is an endomorphism $\phi \colon \mathbb{P}^1 \to \mathbb{P}^1$ satisfying $\phi(x_i) = y_i$ for every $i=1,...,m$.
\end{lemma}
\begin{proof}
Without loss of generality, we may assume that none of the $x_i$ or $y_i$ is the point at infinity. But then a suitable Lagrange interpolation polynomial does the job.
\end{proof}

The basic idea of our proof that $\Sym^m(C \times \mathbb{P}^1)$ is geometrically-special is to take a covering set for $\Sym^m(C)$ and turn it into a covering set for $\Sym^m(C \times \mathbb{P}^1)$ by postcomposing with many different morphisms $\Sym^m(C) \to \Sym^m(C \times \mathbb{P}^1)$ coming from many different morphisms $C \to \mathbb{P}^1$. After these morphisms are constructed, we may test the density of the graphs after projecting down to $\Sym^m(C) \times \Sym^m(\mathbb{P}^1)$, and then it only remains to do the bookkeeping. As this approach does not depend on $C$ being a curve, we state the result in more generality. 

\begin{theorem}\label{thm:sym_n_is_geom_special_0}
Let $m$ be a positive integer and let $X$ be a quasi-projective variety such that $\Sym^m(X)$ is geometrically-special. Then $\Sym^m(X\times \mathbb{P}^1)$ is geometrically-special.
\end{theorem}
\begin{proof}
Since geometric-specialness is a ``birational invariant'' (see \cite[Lemma~2.6]{JR}), we may replace $X$ by a blow-up. Thus, we may assume that there is a dominant morphism $\pi \colon X \to \mathbb{P}^1$. Let $s \in \Sym^m(X)$ be a point through which there is a covering set. We may assume that $s$ represents an $m$-tuple of pairwise distinct points of $X$, say $s = [x_1,...,x_m]$ where we picked an arbitrary ordering. We may even assume that the $\pi(x_k)$ are pairwise distinct (since the set of such points is a non-empty open). We now construct a covering set through the point $[(x_1,z_1),...,(x_m,z_m)] \in \Sym^m(X \times \mathbb{P}^1)$ for any $m$-tuple $(z_1,...,z_m) \in (\mathbb{P}^1)^m$.

Let $(\phi_i \colon (D,d) \to (\Sym^m(X),s))_{i \in I}$ be a covering set. We may shrink the covering set (while retaining its status as a covering set) by removing all morphisms whose image does not contain an $m$-tuple disjoint from the set $\pi^{-1}(\pi(\{x_1,...,x_m\}))$. Now consider the following set:
$$ J := \{(a_1,...,a_m,b_1,...,b_m) \in (\mathbb{P}^1)^{2m}~|~\text{the}~a_k~\text{are pairwise distinct and each distinct from all}~\pi(x_k)\} $$
For every $j = (a_1,...,a_m,b_1,...,b_m) \in J$, we let $\alpha_j$ be any endomorphism of $\mathbb{P}^1$ which sends the points $(a_1,...,a_m,\pi(x_1),...,\pi(x_m))$ to the points $(b_1,...,b_m,z_1,...,z_m)$ (this exists by Lemma \ref{lemma:interpolationpolys}). We obtain morphisms $\beta_j := \alpha_j \circ \pi \colon X \to \mathbb{P}^1$. These give morphisms $(\id_X, \beta_j) \colon X \to X \times \mathbb{P}^1$, which induce morphisms $\gamma_j \colon \Sym^m(X) \to \Sym^m(X \times \mathbb{P}^1)$. Our covering set then consists of the morphisms $(\gamma_j \circ \phi_i)_{(i,j) \in I \times J}$. It remains to verify that this is indeed a covering set.

For this, first note that our base point $d \in D$ always gets mapped to $s = [x_1,...,x_m] \in \Sym^m(X)$ under $\phi_i$. The $x_k$ always get mapped to the corresponding $z_k$ by construction of $\beta_j$. Thus, the image of $d$ in $\Sym^m(X \times \mathbb{P}^1)$ is always $[(x_1,z_1),...,(x_m,z_m)]$, as required.

By Lemma \ref{lemma:finitereflectsdensity2}, the density of the graphs may be tested after projection along $\tau \colon \Sym^m(X \times \mathbb{P}^1) \to \Sym^m(X) \times \Sym^m(\mathbb{P}^1)$. To verify the density now, start by fixing an $i \in I$. By our choice of the covering set for $\Sym^m(X)$, there is a point $d' \in D$ such that $\phi_i(d')$ is an $m$-tuple of pairwise distinct points of $X$ completely disjoint from the set $\pi^{-1}(\pi(\{x_1,...,x_m\}))$. (In fact, since this is an open condition on $d'$, infinitely many such $d'$ exist.) Fixing one $d'$ for now, we see that $\phi_i(d')$ appears, in some ordering, as the first half of an element of $J$. In fact, it does so infinitely many times, as it appears $m!$ times for every tuple $(b_1,...,b_m) \in (\mathbb{P}^1)^m$. This implies:
$$ \overline{\bigcup_{j \in J} (d', (\tau \circ \gamma_j \circ \phi_i)(d'))} = (d', \phi_i(d')) \times \Sym^m(\mathbb{P}^1) \subseteq D \times \Sym^m(X) \times \Sym^m(\mathbb{P}^1) $$
By using either Theorem \ref{thm:geomspeccriterion} or by using that infinitely many such $d'$ exist, we obtain
$$ \overline{\bigcup_{j \in J} \Gamma_{\tau \circ \gamma_j \circ \phi_i}} = \Gamma_{\phi_i} \times \Sym^m(\mathbb{P}^1) \subseteq D \times \Sym^m(X) \times \Sym^m(\mathbb{P}^1)$$
and taking the union over $i \in I$ establishes the required density, since the $\phi_i$ form a covering set for $\Sym^m(X)$.
\end{proof}

\begin{corollary}\label{cor:sym_n_is_geom_special}
Let $C$ be a smooth projective curve of genus $g$ and let $m\geq 1$ be an integer. Then $\Sym^m(C \times \mathbb{P}^1)$ is geometrically-special if and only if $m\geq g$.
\end{corollary}
\begin{proof}  
If $m\geq g$, we have that $\Sym^m(C)$ is geometrically-special by Corollary \ref{cor:sympowergeomspec}, so that $\Sym^m(C\times \mathbb{P}^1)$ is geometrically-special by Theorem\ref{thm:sym_n_is_geom_special_0}.
Conversely, if $\Sym^m(C\times \mathbb{P}^1)$ is geometrically-special, then $\Sym^m(C)$ is geometrically-special (as the former surjects onto the latter) and hence $m \geq g$ by Corollary~\ref{cor:sympowergeomspec}.
\end{proof}

Let $E$ be an elliptic curve. To prove the geometric-specialness of $\Sym^m(C \times E)$ we will use that the existence of a nonconstant morphism $C \to E$ implies that $\Sym^m(C \times E)\to \Sym^m(C)$ has many sections. We will then postcompose the covering sets through well-chosen points of $\Sym^m(C)$ with these sections to obtain covering sets for $\Sym^m(C \times E)$. 

\begin{theorem}\label{thm:symCcrossE_geom_special}
Let $C$ be a smooth projective curve of genus $g$ and let $E$ be an elliptic curve admitting a surjection $\pi \colon C \to E$. Let $m \geq g$ be a natural number. Then $\Sym^m(C \times E)$ is geometrically-special.
\end{theorem}
\begin{proof}
Let $[c_1,...,c_m] \in \Sym^m(C)$ be a point such that $\pi(c_i) \in E$ is a torsion point for every $i=1,...,m$ and such that there is a covering set for $\Sym^m(C)$ through $[c_1,...,c_m]$. Observe that the set of such points is dense in $\Sym^m(C)$ as the first condition holds on a dense set and the second condition holds on a nonempty open by Corollary \ref{cor:sympowergeomspec}. Fix an integer $k$ and let $z_i = [k](c_i)$, where $[k] \colon E \to E$ denotes the multiplication-by-$k$ morphism. Note that the set of all points $[(c_1,z_1),...,(c_m,z_m)] \in \Sym^m(C \times E)$ obtained this way is dense in $\Sym^m(C \times E)$. Thus, to show geometric specialness, it suffices to construct a covering set for $\Sym^m(C \times E)$ through such a point $[(c_1,z_1),...,(c_m,z_m)] \in \Sym^m(C \times E)$. Since we assumed the $\pi(c_i)$ to be torsion points of $E$, there is an integer $n$ such that $[n](\pi(c_i)) = 0 \in E$ for all $i=1,...,m$.

By construction, there is a covering set $(\phi_i \colon (D,d) \to (\Sym^m(C), [c_1,...,c_m]))_{i \in I}$. For each integer $j \in \mathbb{Z}$, we define the morphism $\gamma_j \colon \Sym^m(C) \to \Sym^m(C \times E)$ to be the $m$-th symmetric power of the morphism $(\id_C,[nj+k] \circ \pi) \colon C \to C \times E$. We claim that the family of morphisms $(\gamma_j \circ \phi_i \colon D \to \Sym^m(C \times E))_{i \in I, j \in \mathbb{Z}}$ constitutes a covering set for $\Sym^m(C \times E)$ through the point $[(c_1,z_1),...,(c_m,z_m)]$.

To verify this, first note that $\gamma_j$ sends the point $[c_1,...,c_m]$ to $[(c_1,[nj+k](\pi(c_1))),...,(c_m,[nj+k](\pi(c_m)))]$, and as we assumed all $\pi(c_i)$ to be $n$-torsion, we have $[nj+k](\pi(c_i)) = [k](\pi(c_i)) = z_i$. As $\phi_i$ sends the point $d \in D$ to $[c_1,...,c_m] \in \Sym^n(C)$ by definition, this implies that the morphisms $\gamma_j \circ \phi_i$ do indeed send $d \in D$ to $[(c_1,z_1),...,(c_m,z_m)] \in \Sym^n(C \times E)$. It remains to verify the density of the graphs in $D \times \Sym^m(C \times E)$. 

Next, we verify that the morphisms $(\gamma_j \colon \Sym^m(C) \to \Sym^m(C \times E))_{j \in \mathbb{Z}}$ have jointly dense image. To see this, first note that by Lemma \ref{lemma:finitereflectsdensity2}, we may test this after projecting to $\Sym^m(C) \times \Sym^m(E)$. Next, note that by Theorem \ref{thm:geomspeccriterion}, it suffices to show that there is a point $x \in \Sym^m(C)$ such that the set $\{ \Sym^m([nj+k] \circ \pi)(x)~|~j \in \mathbb{Z} \}$ is dense in $\Sym^m(E)$. To see that such a point $x$ exists, let $e \in E^m$ be a nondegenerate point and choose $x \in \Sym^m(C)$ such that $\Sym^m(\pi)(x)$ is the image of $e$ in $\Sym^m(E)$. This $x$ then has the desired property.

To conclude, observe that
\[ \bigcup_{i \in I, j \in \mathbb{Z}} \Gamma_{\gamma_j \circ \phi_i} = \bigcup_{j \in \mathbb{Z}} (\id_D, \gamma_j) \left(\bigcup_{i \in I} \Gamma_{\phi_i}\right) \]
so that
\[ \overline{\bigcup_{i \in I, j \in \mathbb{Z}} \Gamma_{\gamma_j \circ \phi_i}} = \overline{\bigcup_{j \in \mathbb{Z}} \overline{(\id_D, \gamma_j) \left(\bigcup_{i \in I} \Gamma_{\phi_i}\right)}} \]
As the $\phi_i$ form a covering set, it follows that 
\[ \overline{\bigcup_{i \in I, j \in \mathbb{Z}} \Gamma_{\gamma_j \circ \phi_i}} = \overline{\bigcup_{j \in \mathbb{Z}} (\id_D, \gamma_j)(D \times \Sym^m(C))} = D \times \overline{\bigcup_{j \in \mathbb{Z}} \gamma_j(\Sym^m(C))} \]
and since we verified that the morphisms $\gamma_j$ have jointly dense image, we conclude.
\end{proof}

\section{Potential density} \label{section:potential_density}

In this section we  first characterize which symmetric powers of $C\times \mathbb{P}^1$   are arithmetically-special (i.e., have a potentially dense set of rational points). In our approach, we will need the existence of rational points on certain twists of $(\mathbb{P}^1)^m$. This naturally leads us to studying $L$-rational points on $(\mathbb{P}^1)^m$ whose coordinates form a transitive $S_m$-set. 

\begin{lemma} \label{lem:density}
Let $K$ be an infinite field and let $L$ be a finite separable field extension of $K$ of degree $m$. Let $\tau_1,\ldots, \tau_m \colon L\to \overline{K}$ be the $m$ pairwise distinct embeddings of $L$ into $\overline{K}$. For $\alpha$ in $L$, let $P_{\alpha} = (\tau_1(\alpha),\ldots, \tau_m(\alpha)) \in \mathbb{A}^m(\overline{K})$. Then the set 
\[ R = \{P_\alpha~|~\alpha\in L \} \subset \mathbb{A}^m(\overline{K}) \]
is dense in $\mathbb{A}^m$.
\end{lemma}
\begin{proof}
Let $\alpha_1,\ldots,\alpha_m$ be a $K$-basis for $L$. Consider the $(m \times m)$-matrix $M = (\tau_i\alpha_j)$. As is well-known, $(\det M)^2$ is the discriminant of the $K$-basis $\alpha_1,\ldots, \alpha_m$, which is nonzero since $L/K$ is separable. Then the matrix $M$ defines an invertible linear map $\mathbb{A}^m \to \mathbb{A}^m$ under which $R$ is the image of $\mathbb{A}^m(K)$. Since $K$ is infinite, $\mathbb{A}^m(K)$ is dense in $\mathbb{A}^m$, and it follows that $R$ is also Zariski-dense in $\mathbb{A}^m$.
\end{proof}

Note that Lemma \ref{lem:density} gives a simple proof of the Primitive Element Theorem in the case of infinite fields. Indeed, the set associated to non-primitive elements $\{ P_\alpha~|~\alpha\in L, L\neq K(\alpha)\}$ is not dense in $\mathbb{A}^m$ (it's contained in the union of hyperplanes of the form $x_i=x_j$, $i\neq j$). Therefore there must exist a primitive element for $L/K$. Note that this proof depends only on the following two facts: $\mathbb{A}^m(K)$ is dense in $\mathbb{A}^m$ if $K$ is infinite, and the discriminant of any $K$-basis of a finite separable extension $L/K$ is nonzero.

\begin{proposition}\label{prop:S_n_improved}
Let $K$ be an infinite field and let $K \subseteq L$ be a finite Galois extension. Choose an embedding $\Gal(L/K) \subseteq S_m$ for some integer $m$. Then the set
\[ \{ (x_1,\ldots,x_m) \in (\mathbb{P}^1)^m(L)~|~\text{for all } i=1,\ldots,m \text{ and all } \sigma \in \Gal(L/K), \text{ we have } \sigma(x_i)=x_{\sigma(i)} \} \]
is dense in $(\mathbb{P}^1)^m$.
\end{proposition}
\begin{proof}
Let $G$ be the image of $\Gal(L/K)$ in $S_m$. We first treat the case that $G$ is a transitive subgroup of $S_m$. In this case, let $G' = S_{m-1} \cap G$, where we embed $S_{m-1} \subseteq S_m$ as the stabilizer of a point, and let $K' = L^{G'}$ be the corresponding fixed field. Then $K \subseteq K'$ is an extension of degree $m$ with Galois closure $L$. Let $\tau_1,\ldots,\tau_m$ be the $m$ distinct embeddings of $K'$ in $L$ over $K$. Then $G$ acts on the set $\{\tau_1,\ldots,\tau_m\}$ and after renumbering the $\tau_i$, we may assume that $\sigma \circ \tau_i = \tau_{\sigma(i)}$ for every $i=1,\ldots,m$ and every $\sigma \in G$. Then, by Lemma \ref{lem:density}, we have that
\[ \{ (\tau_1(\alpha),\ldots,\tau_m(\alpha))~|~\alpha \in K' \} \]
is a dense set of elements of $(\mathbb{P}^1)^m(L)$ with the desired transformation behaviour under $\Gal(L/K)$, finishing the proof if $G \subseteq S_m$ is transitive. 

If $G \subseteq S_m$ is not a transitive subgroup, let $r_1,\ldots,r_l$ denote the sizes of the orbits. After renumbering, we may assume that the orbits are $\{1,\ldots,r_1\}, \{r_1+1,\ldots,r_1+r_2\},$ and so on. For $j=1,\ldots,l$, let $G_j$ be the image of $G$ under the natural restriction homomorphism $G \to S_{r_j}$ and let $N_j \subseteq G$ be the kernel of $G \to G_j$. Let $L_j \subseteq L$ be the fixed field of $N_j$. Then $K \subseteq L_j$ is a Galois extension with Galois group $G_j \subseteq S_{r_j}$ and the subgroup $G_j \subseteq S_{r_j}$ is transitive. Thus, by the first paragraph of this proof, the following set is dense in $(\mathbb{P}^1)^{r_j}$.
\[ \Sigma_j := \{ (x_1,...,x_{r_j}) \in (\mathbb{P}^1)^{r_j}(L_j)~|~\text{for all } i=1,\ldots,r_j \text{ and all } \sigma \in G_j, \text{ we have } \sigma(x_i)=x_{\sigma(i)} \} \]
Thus, the product set $\Sigma := \Sigma_1 \times \cdots \times \Sigma_l \subseteq (\mathbb{P}^1)^{r_1}(L_1) \times \cdots \times (\mathbb{P}^1)^{r_l}(L_l) \subseteq (\mathbb{P}^1)^m(L)$ is dense in $(\mathbb{P}^1)^m$. Now note that by construction, the elements of $\Sigma$ have the desired transformation behaviour under $\Gal(L/K)$, finishing the proof in general. 
\end{proof}

The previous proposition will provide an elementary proof of the density of $K$-points on certain twists of $(\mathbb{P}^1)^m$ appearing in our proof of Theorem \ref{thm:sym_X_times_P1} below. We will see later that this density can also be proven using the structure of such twists as twisted flag varieties; see the proof of Lemma \ref{lem:WHP_A_points_on_Sym}.

If $X\to S$ is a quasi-projective morphism of noetherian schemes and $m\geq 1$ is an integer, then $S_m$ acts on the fiber product $X^m = X \times_S \ldots \times_S X$. We will denote its quotient by $\Sym^m_S(X)$; note that this is again a quasi-projective scheme over $S$. This follows from \cite[Theorem~V.4.1]{SGA3}, as explained in \cite[Remarque~V.5.1]{SGA3}.

\begin{lemma}\label{lem:A_points_on_Sym}
Let $A$ be a noetherian integral domain with infinite fraction field $K$, and let $\mathcal{X}$ be a quasi-projective integral scheme over $A$. If $m \geq 1$ is any integer for which $\Sigma := \Sym^m_A(\mathcal{X})(A)^{(1)}$ is dense in $\mathcal{X}$, then $\Sym^m_A(\mathcal{X}\times \mathbb{P}^1_A)(A)^{(1)}$ is dense in $\Sym^m_A(\mathcal{X}\times \mathbb{P}^1_A)$.
\end{lemma}
\begin{proof}
Let $X=\mathcal{X}_K$ be the generic fiber of $\mathcal{X}\to \Spec A$. Note that $\mathcal{X}^m\to \Sym^m_A(\mathcal{X})$ extends the natural quotient morphism $X^m \to \Sym^m(X)$ over $A$. 

For every subgroup $H \subseteq S_m$, let $\Sigma_H$ be the set of points $t \in \Sigma$ for which the fiber $(X^m)_t$ is reduced and every connected component of the finite $K$-scheme $(X^m)_t$ is an $H$-torsor. Note that every point of $\Sigma$ not lying on the big diagonal lies in one of the $\Sigma_H$. Thus, as we assumed $\Sigma$ to be dense, we see that $\Sym^m(X)$ is equal to the (finite) union of the closures $\overline{\Sigma_H}$ and the big diagonal. Consequently, since $\Sym^m(X)$ is irreducible (and since the big diagonal is not dense), we conclude that $\Sym^m(X)= \overline{\Sigma_H}$ for some subgroup $H\subseteq S_m$. In other words, there is a subgroup $H \subseteq S_m$ for which $\Sigma_H \subseteq \Sym^m(X)$ is dense.

For every $t \in \Sigma_H$, we claim that the fiber $F_t$ of the natural projection $\Sym^m_A(\mathcal{X} \times \mathbb{P}^1_A)\to \Sym^m_A(\mathcal{X})$ over $t \in \Sym^m_A( \mathcal{X})(A)^{(1)}$ has a dense set of near-integral $A$-points. To show this, it suffices to show that $F_t(K)$ is dense, as $F_t$ is a proper $A$-scheme. To do so, fix a $t \in \Sigma_H$ and write $L$ for the field extension given by the connected components of $(X^m)_t$. The field extension $L/K$ is a finite Galois extension. Observe that $F_t \otimes_K L$ is isomorphic to $(\mathbb{P}^1_L)^m$ and that the image of the set
\[ \{ (x_1,\ldots,x_m) \in (\mathbb{P}^1)^m(L)~|~\text{for all } i=1,\ldots,m \text{ and all } \sigma \in \Gal(L/K), \text{ we have } \sigma(x_i)=x_{\sigma(i)} \} \]
in $F_t(L)$ is contained in $F_t(K)$. Hence it follows from Proposition \ref{prop:S_n_improved} that $F_t(K)$ is dense.

We thus have shown that there is a dense set of $t \in \Sym^m_A(\mathcal{X})(A)^{(1)}$ for which the fiber of $\Sym^m_A(\mathcal{X} \times \mathbb{P}^1_A) \to \Sym^m_A(\mathcal{X})$ over $t$ has a dense set of near-integral $A$-points. We conclude that $\Sym^m_A(\mathcal{X}\times \mathbb{P}^1_A)(A)^{(1)}$ is dense, as required.
\end{proof}

\begin{theorem}\label{thm:sym_X_times_P1}
Let $k$ be an algebraically closed field of characteristic zero. Let $X$ be a quasi-projective variety over $k$ such that $\Sym^m(X)$ is arithmetically-special over $k$. Then $\Sym^m(X\times \mathbb{P}^1)$ is arithmetically-special over $k$.
\end{theorem}
\begin{proof}
Let $A\subset k$ be a $\mathbb{Z}$-finitely generated subring with fraction field $K$, and let $\mathcal{X}$ be a quasi-projective model for $X$ over $A$. Replacing $\Spec A$ by a dense affine open if necessary, we may assume that $\Sym^m_A(\mathcal{X})$ is a quasi-projective model for $\Sym^m(X)$ over $A$ (or, alternatively, we can avoid spreading out by simply appealing to the aforementioned result in \cite{SGA3}). Since $\Sym^m(X)$ is arithmetically-special over $k$, replacing $A$ by a suitable finitely generated extension, we may assume that $\Sigma := \Sym_A^m(\mathcal{X})(A)^{(1)} \subseteq \Sym^m(X)(k)$ is dense. It now follows from Lemma \ref{lem:A_points_on_Sym} that $\Sym^m_A(\mathcal{X}\times \mathbb{P}^1_A)(A)^{(1)}$ is dense. 
\end{proof}

\begin{lemma}\label{lemma:sym_m_curve_is_arithmetically_special}
Let $k$ be an algebraically closed field of characteristic zero. Let $m$ be a positive integer and let $C$ be a smooth projective curve of genus $g$ over $k$. Then $\Sym^m(C)$ is arithmetically-special over $k$ if and only if $m\geq g$. 
\end{lemma}
\begin{proof}
If $m < g$, then the image of $\Sym^m(C)\to \Pic^m(C)$ is a positive-dimensional closed subvariety of an abelian variety of general type, and thus not arithmetically-special by Faltings's theorem \cite{FaltingsLang}. It follows that $\Sym^m(C)$ is not arithmetically-special. If $m\geq g$, note that $\Sym^m(C)$ is birational to $\mathbb{P}^{m-g}\times \Pic^m(C)$ (see Lemma \ref{lemma:sym_n_coherent_sheaf}). Since $\mathbb{P}^{m-g}$ and $\Pic^m(C)$ are arithmetically-special, so is their product $\mathbb{P}^{m-g}\times \Pic^m(C)$. Since being arithmetically-special is a birational invariant, we conclude that $\Sym^m(C)$ is arithmetically-special.
\end{proof}

\begin{corollary} \label{cor:sym_n_is_arith_special}
Let $k$ be an algebraically closed field of characteristic zero. Let $m$ be an integer and let $C$ be a smooth projective curve of genus $g$ over $k$. Then $\Sym^m(C\times \mathbb{P}^1_k)$ is arithmetically-special if and only if $m \geq g$. 
\end{corollary}
\begin{proof}
Note that the natural projection $C \times \mathbb{P}^1_k \to C$ induces a surjective morphism $\Sym^m(C\times \mathbb{P}^1_k)\to \Sym^m(C)$. In particular, if $\Sym^m(C\times \mathbb{P}^1_k)$ is arithmetically-special, then $\Sym^m(C)$ is arithmetically-special, so that $m\geq g$ by Lemma \ref{lemma:sym_m_curve_is_arithmetically_special}. Conversely, if $m\geq g$, then $\Sym^m(C)$ is arithmetically-special by Lemma \ref{lemma:sym_m_curve_is_arithmetically_special}, so that the result follows from Theorem \ref{thm:sym_X_times_P1}. 
\end{proof}

We can now show that $\Sym^m(C \times \mathbb{P}^1)$ for $C$ a smooth projective curve of genus $g \geq 2$ and $m \geq g$ also provides a counterexample to Hassett--Tschinkel's arithmetic puncturing problem (Problem~\ref{pAP}). That is, we can now prove Theorem~\ref{thm:ce_ar} from the introduction.

\begin{proof}[Proof of Theorem~\ref{thm:ce_ar}]
The variety $\Sym^m(C \times \mathbb{P}^1)$ is arithmetically-special by Corollary~\ref{cor:sym_n_is_arith_special}; this shows $(1)$.
The complement of the big diagonal in $\Sym^m(C \times \mathbb{P}^1)$ is not arithmetically-special by Theorem~\ref{thm:sym_prods_2}; this shows $(2)$.
That $\Sym^m(C \times \mathbb{P}^1)$ has canonical singularities was already shown in Remark~\ref{thm:ce_geom} (see Section~\ref{section:ascending_descending} for the proof).
\end{proof}

We now prove the potential density of rational points on $\Sym^m(C\times E)$ when $m$ is at least the genus of $C$, assuming $E$ is an elliptic curve and $C$ admits a cover $C\to E$.

\begin{theorem} \label{thm:symCcrossE_arith_special}
Let $E$ be an elliptic curve over a finitely generated field $K$ of characteristic zero and let $C$ be a smooth projective curve of genus $g$ over $K$. Assume that $C_{\overline{K}}$ dominates $E_{\overline{K}}$. If $m\geq g$ is a positive integer, then there is a finite field extension $L/K$ such that  $\Sym^m(C \times E)(L)$ is dense in $\Sym^m(C\times E)$, 
\end{theorem}
\begin{proof} 
Replacing $K$ by a finite field extension if necessary, we may assume that there is a surjective morphism $\pi\colon C\to E$. The morphism $\pi$ then induces a natural morphism $\id_C \times \pi \colon C \to C \times E$. Let $[n]$ be the self-map of $C\times E$ given by multiplication with $n$ on $E$ and the identity on $C$, i.e., $[n]\colon C\times E\to C\times E$ sends $(c,x)$ to $(c,nx)$. By the functoriality of symmetric products, the composed morphism $[n] \circ (\id_C \times \pi) \colon C\to C\times E$ induces a morphism $\pi_n \colon \Sym^m(C)\to \Sym^m(C\times E)$. 
We let $\phi_n\colon \Sym^m(C)\to \Sym^m(E)$ be the morphism $\pi_n$ composed with the projection $\Sym^m(C\times E)\to \Sym^m(E)$. Note that $\phi_n$ is also the morphism induced by $[n] \circ \pi \colon C \to E$. Since $E^m$ is an abelian variety, there is a point $g\in E^m(\overline{K})$ such that the subgroup generated by $g$ in $E^m$ is dense. (Such a point is called a non-degenerate point of $E^m$.) It follows that, replacing $K$ by a finite field extension if necessary, there is a point $y_0$ in $\Sym^m(C)(K)$ such that the set $\{\phi_i(y_0) \ | \ i=1,2,\ldots\} $ is dense. (Take $y_0$ to be any point mapping to the class of $g$ via $\phi_1$.) We have thus the following commutative diagram:

\[
\xymatrix{ 
 \Sym^m(C\times E)\ar[d] \ar[rr]^{\mathrm{projection}} & & \Sym^m(E) \\ 
 \Sym^m(C)\times \Sym^m(E) \ar[d] & & \\ 
 \Sym^m(C) \ar@/{}^{5pc}/^{\pi_n}[uu] \ar[uurr]_{\phi_n} & &
}
\]

Replacing $K$ by a finite field extension if necessary, we may assume that $\Sym^m(C)$ has a dense set of $K$-points. Let $Z_n$ be the image of $\pi_n$, and note that $Z_n$ is isomorphic to $\Sym^m(C)$. In particular, the set of $K$-points $Z_n(K)$ is dense in $Z_n$. We claim that $\cup_{n=1}^\infty Z_n(K)$ is dense in $\Sym^m(C\times E$). To prove this, it suffices to show that its image in $\Sym^m(C) \times \Sym^m(E)$ is dense. 
To verify this, define $Y := \Sym^m(C)$ and $X := \Sym^m(E)$. Since $Y(K)$ is dense in $Y$ and $\{\phi_n(y_0)~|~n=1,2,\ldots\}$ is dense in $X$, by Corollary \ref{cor:geomspeccriterion_arithmetic}, the subset $\cup_{n=1}^\infty \Gamma_{\phi_n}(K)$ is dense in $Y \times X = \Sym^m(C)\times \Sym^m(E)$.
\end{proof}

\section{The Hilbert property}

Recall that a proper variety $X$ over a field $k$ is said to have the \emph{Hilbert property over $k$} if $X(k)$ is not thin \cite[\S 3]{SerreTopicsGalois}. Concretely, we have that $X$ has the Hilbert property over $k$ if, for every finite collection of finite surjective morphisms $(\pi_i \colon Y_i \to X)_{i=1}^n$ with $Y_i$ a normal (integral) variety over $k$ and $\deg \pi_i \geq 2$, the set $X(k)\setminus \cup_{i=1}^n \pi_i(Y_i(k))$ is dense in $X$.
 
Recall that a field $K$ is \emph{Hilbertian} if $\mathbb{P}^1_K$ has the Hilbert property over $K$. For example, every number field is Hilbertian \cite[\S3.4]{SerreTopicsGalois}. 
We will use that a twist of $(\mathbb{P}^1_K)^m$ satisfies the Hilbert property if it has a $K$-point; this follows from Bary-Soroker--Fehm--Petersen's result that any smooth compactification of a linear algebraic group over a number field has the Hilbert property \cite[Corollary~4.2]{BSFP}. 

\begin{theorem}\label{thm:hp_for_twists_of_powers_of_P1}
Let $K$ be a Hilbertian field of characteristic zero and let $X$ be a smooth proper variety such that $X_{\overline{K}}$ is isomorphic to a power of $\mathbb{P}^1_{\overline{K}}$. If $X(K)$ is non-empty, then $X$ has the Hilbert property over $K$.
\end{theorem}
\begin{proof}
Note that $X$ is a \emph{twisted flag variety} over $K$ (in the sense of \cite[\S 6, Definition~1]{Demazure}).
Let $x\in X(K)$ be a $K$-rational point. Let $G=\underline{\Aut}_{X/k}$ be the automorphism group scheme of $X$. By \cite[Proposition~4]{Demazure}, we have that $G$ is a smooth affine finite type group scheme over $K$ whose connected component $G^0$ is a connected semisimple linear algebraic group over $K$ such that $X$ is a homogeneous space under $G$. 
Let $P$ be the stabilizer group scheme of $x$. Then $P$ is a parabolic subgroup, hence connected. Thus, as $K$ is a perfect Hilbertian field, the homogeneous space $G/P$ has the Hilbert property over $K$ by \cite[Corollary~4.6]{BSFP}. Since $X$ is isomorphic to $G/P$ (via $G\to X$ defined by $g\mapsto g\cdot x$), we conclude that $X$ has the Hilbert property over $K$. 
\end{proof}

A smooth projective variety over $k$ with the Hilbert property over $k$ has a dense set of $k$-points. However, the converse fails. For example, if $E$ is an elliptic curve of positive rank over a number field $k$, then $E$ does not have the Hilbert property over $k$ (despite $E(k)$ being dense). It does however satisfy the weak Hilbert property (by Faltings's theorem \cite{Faltings2-2}).
In fact, a smooth proper variety over a number field $k$ has the Hilbert property if and only if it has the weak Hilbert property and $X_{\overline{k}}$ has no non-trivial finite \'etale covers (see \cite{CZHP}).

The weak Hilbert property for $X$ guarantees that given a ramified cover $Y \to X$, many fibers $Y_x$ do not have a $k$-point. Assuming that $Y \to X$ is Galois and ``genuinely ramified'', this statement can be strengthened as follows.

\begin{lemma}\label{lemma:Galois_genuinely_ramified}
Let $X$ be a smooth proper variety with the weak Hilbert property over $K$. Let $Y\to X$ be a ramified Galois covering which has no nontrivial \'etale subcovering. Then, the set of $x\in X(K)$ such that $Y_x$ is integral is dense.
\end{lemma}
\begin{proof}
Let $G$ be the Galois group of $Y\to X$. Consider the collection of coverings $Y/H \to X$ as $H$ runs over all subgroups $H\neq G$. Note that each such covering is ramified (as it is a subcovering of $Y\to X$). Therefore, by applying the weak Hilbert property to the collection $(Y/H\to X)_{H \subset G, H\neq G}$, we see that the set of non-branch points $x$ in $X(K)$ with $(Y/H)_x(K)=\emptyset$ for every $H\subsetneq G$ is dense. Note that for each such $x$, the fiber $Y_x$ is integral.
\end{proof} 

We will use that the weak Hilbert property is inherited by the total space of a family of varieties satisfying the Hilbert property over a base satisfying the weak Hilbert property. The precise result we need is a consequence of a general fibration theorem proven in \cite{Luger3} (which improves on the fibration theorems of \cite{BSFP} and \cite{Jnef}).

\begin{theorem}[Mixed fibration theorem] \label{thm:mixed_fibration_theorem}
Let $K$ be a field of characteristic zero and let $f \colon X \to S$ be a proper surjective morphism of normal varieties over $K$. Let $\Gamma \subset X(K)$ be a subset. Let $\Sigma \subset S(K)$ be a subset which is not strongly thin. Suppose that, for every $s$ in $\Sigma$, the proper $K$-scheme $X_s$ is integral and normal and that the subset $X_s(K) \cap \Gamma$ is not thin in $X_s$. Then $\Gamma$ is not strongly thin in $X$.
\end{theorem}

\subsection{Symmetric products} 

\begin{lemma}\label{lemma:whp} 
Let $C$ be a smooth projective curve of genus $g$ over a finitely generated field $K$ of characteristic zero. Then there is a finite field extension $L/K$ such that the smooth projective variety $\Sym^m(C_L)$ has the weak Hilbert property over $L$ if and only if $m\geq g$. 
\end{lemma}
\begin{proof}
If there is a finite field extension $L/K$ such that $\Sym^m(C_L)$ has the weak Hilbert property over $L$, then $\Sym^m(C_{\overline{K}})$ is arithmetically-special (trivially), so that $m\geq g$ by Lemma \ref{lemma:sym_m_curve_is_arithmetically_special}. Now, assume $m\geq g$, and let $\mathrm{Jac}(C)$ be the Jacobian of $C$. Since $m\geq g$, replacing $K$ by a finite field extension if necessary, we may assume that $\Sym^m(C)$ is $K$-birational to $\mathbb{P}^{m-g}_K\times \mathrm{Jac}(C)$ (by Lemma \ref{lemma:sym_n_coherent_sheaf}).
Note that $\mathbb{P}^{m-g}$ has the Hilbert property over $K$ \cite[\S3]{SerreTopicsGalois}. Moreover, replacing $K$ by a finite field extension if necessary, by work of Frey--Jarden \cite{FreyJarden}, the abelian variety $\mathrm{Jac}(C)$ has a dense set of $K$-points (see \cite[Corollary~3.8]{JAut} for a precise statement), and thus the weak Hilbert property over $K$ \cite{CDJLZ}. In particular, by the mixed fibration theorem (Theorem \ref{thm:mixed_fibration_theorem}) (or the product theorem for WHP \cite[Theorem~1.9]{CDJLZ}), the variety $\mathbb{P}^{m-g}_K\times \mathrm{Jac}(C)$ has the weak Hilbert property over $K$. In particular, since the weak Hilbert property is a birational invariant amongst smooth projective varieties \cite[Proposition~3.1]{CDJLZ}, it follows that $\Sym^m(C)$ has the weak Hilbert property over $K$, as required.
\end{proof}

As an interesting application of the weak Hilbert property of $\Sym^m(C)$, we obtain the infinitude of $S_m$-points on curves for every $m$ at least the genus:

\begin{corollary}\label{cor:Sm_points_on_C}
Let $C$ be a smooth projective geometrically connected curve over a finitely generated field $K$ of characteristic zero. If $m\geq g$, then there is a finite field extension $L/K$ such that the set of $c$ in $C_L$ whose residue field $\kappa(c)$ is an $S_m$-Galois extension of $L$ is infinite.
\end{corollary}
\begin{proof}
Replacing $K$ by a finite field extension if necessary, we may assume that $\Sym^m(C)$ has the weak Hilbert property over $K$ (Lemma \ref{lemma:whp}). Now, note that the morphism $C^m \to \Sym^m(C)$ has no non-trivial \'etale subcovers. Indeed, for every $x$ in $C(\overline{K})$, the fiber over $[(x, x,\ldots,x)]$ is the single point $(x,\ldots,x)$. In particular, since the morphism $C^m\to \Sym^m(C)$ is generically an $S_m$-torsor, the corollary follows from Lemma \ref{lemma:Galois_genuinely_ramified}.
\end{proof}

\begin{remark}[Wittenberg]\label{remark:wittenberg}
If $m\geq 2g$, then one can prove Corollary \ref{cor:Sm_points_on_C} without appealing to the weak Hilbert property of abelian varieties \cite{CDJLZ}.  We thank Olivier Wittenberg for allowing us to include the following argument.

First, extending $K$ if necessary, we may assume that $\Pic^m(C)(K)$ is dense.
Let $p$ be a general point of $\Pic^m(C)$. Then the fiber $\Sym^m(C)_p$ of $\Sym^m(C) \to \Pic^m(C)$ over $p$ is a projective space and the cover $(C^m)_p \to (\Sym^m(C))_p$ is generically an $S_m$-torsor. Now consider the cover $C \times \Sym^{m-1}(C) \to \Sym^m(C)$, which is an intermediate cover of $C^m \to \Sym^m(C)$. Note that, passing to the fiber over $p \in \Pic^m(C)$, the projection onto the first coordinate $(C \times \Sym^{m-1}(C))_p \to C$ is a projective bundle with fibers of dimension $m-g$. In particular, we see that $(C \times \Sym^{m-1}(C))_p$ is geometrically irreducible. Furthermore, the covering $C \times \Sym^{m-1}(C) \to \Sym^m(C)$ ramifies only over the big diagonal and the fiber over the generic point of the big diagonal has one point of multiplicity two and is otherwise étale. Consequently, the local monodromy is generated by a single transposition. As $p$ was general, we see that the same holds around the codimension one points of the branch locus of $(C \times \Sym^{m-1}(C))_p \to (\Sym^m(C))_p$. Now, as $(\Sym^m(C))_p$ is geometrically simply connected and $(C \times \Sym^{m-1}(C))_p$ is geometrically irreducible, we see that the local monodromy groups generate the global monodromy group. As a transitive subgroup of $S_m$ generated by transpositions must be the entire symmetric group $S_m$, we see that the global monodromy group of $(C \times \Sym^{m-1}(C))_p \to (\Sym^m(C))_p$ is given by $S_m$. Consequently, its Galois closure is $(C^m)_p \to (\Sym^m(C))_p$ and $(C^m)_p$ is geometrically irreducible. In particular, by Hilbert's irreducibility theorem applied to the projective space $(\Sym^m(C))_p$, we see that $(C^m)_p$ has a dense set of closed points whose Galois group is $S_m$. As $p$ was a general point, the same follows for $C^m$ and hence for $C$. 
\end{remark}

\begin{lemma}\label{lem:WHP_A_points_on_Sym} 
Let $A$ be a noetherian integral domain whose fraction field $K$ is Hilbertian and of characteristic zero. Let $\mathcal{X}$ be a quasi-projective integral scheme over $A$. If $\Sym^m_A(\mathcal{X})(A)^{(1)} \subset \Sym^m(\mathcal{X}_K)(K)$ is not strongly thin, then $\Sym^m_A(\mathcal{X}\times_A \mathbb{P}^1_A)(A)^{(1)}$ is not strongly thin.
\end{lemma} 
\begin{proof}
(We adapt the proof of Theorem \ref{thm:sym_X_times_P1}.)
Define $X = \mathcal{X}_K$ and $\Gamma := \Sym^m_A(\mathcal{X} \times_A \mathbb{P}^1_A)(A)^{(1)}$, and consider the proper surjective morphism $\Sym^m_A(\mathcal{X}\times_A \mathbb{P}^1_A)\to \Sym^m_A(\mathcal{X})$. Choose a dense open $\mathcal{U}\subseteq \mathcal{X}$ and a dominant morphism $\mathcal{U}\to \mathbb{P}^1_A$. This induces a section of $\Sym^m_A(\mathcal{X}\times_A\mathbb{P}^1_A)\to \Sym^m_A(\mathcal{X})$ over $\Sym^m(\mathcal{U})$. 
Define $\Sigma := \Sym^m_A(\mathcal{X})(A)^{(1)} \cap \Sym^m_A(\mathcal{U})(A)^{(1)}$, and note that $\Sigma$ is not strongly thin (as $\Sym^m_A(\mathcal{X})(A)^{(1)}$ is not strongly thin). Now, for every $s$ in $\Sigma$, the fiber of the morphism $\Sym^m(X\times \mathbb{P}^1_K)\to \Sym^m(X)$ over $s$ (where we view $s$ as a $K$-point of $\Sym^m(X)$) is a twist of $(\mathbb{P}^1_K)^m$ with a $K$-point (since $\Sym^m_A(\mathcal{X}\times_A\mathbb{P}^1_A)\to \Sym^m_A(\mathcal{X})$ has a section over $\Sym^m(\mathcal{U})$). Therefore, for such an $s$, the fiber $X_s$ has the Hilbert property over $K$ by Theorem \ref{thm:hp_for_twists_of_powers_of_P1}, i.e., $X_s(K)$ is not thin. Now since $X_s(K)\cap \Gamma = X_s(K)$ by definition of $\Gamma$, we see that $X_s(K)\cap \Gamma$ is not thin. Therefore, by the mixed fibration theorem (Theorem \ref{thm:mixed_fibration_theorem}), we conclude that $\Gamma$ is not strongly thin, as required.
\end{proof}
 
For the sake of clarity, we state the following consequence of Lemma \ref{lem:WHP_A_points_on_Sym}.
\begin{corollary}\label{corollary:sym_x_sym_x_cross_p1}
Let $k$ be an algebraically closed field of characteristic zero. Let $X$ be a quasi-projective normal variety over $k$ such that $\Sym^m(X)$ has the arithmetic weak Hilbert property over $k$. Then $\Sym^m(X\times \mathbb{P}^1_k)$ has the arithmetic weak Hilbert property over $k$.
\end{corollary}
\begin{proof}
Choose suitable models and apply Lemma \ref{lem:WHP_A_points_on_Sym}.
\end{proof}

\begin{theorem} \label{thm:whp_for_sym}
Let $K$ be a finitely generated field of characteristic zero. Let $m \geq g$ be an integer and let $C$ be a smooth projective curve of genus $g$ over $K$. Then there is a finite field extension $L/K$ such that $\Sym^m(C_L\times \mathbb{P}^1_L)$ has the weak Hilbert property over $L$.
\end{theorem}
\begin{proof}
Since $m\geq g$, replacing $K$ by a finite field extension if necessary, we may and do assume that $\Sym^m(C)$ has the weak Hilbert property (Lemma \ref{lemma:whp}). It now follows from Corollary \ref{corollary:sym_x_sym_x_cross_p1} that, replacing $K$ by a finite field extension if necessary, the normal variety $\Sym^m(C\times \mathbb{P}^1_K)$ has the weak Hilbert property over $K$.
\end{proof}

We note that the prediction made by Conjecture \ref{conj:whp} is that some resolution of singularities of $\Sym^n(C\times \mathbb{P}^1_K)$ has the weak Hilbert property. This follows however from the fact that $\Sym^n(C\times \mathbb{P}^1_K)$ has the weak Hilbert property and the following lemma.

\begin{lemma}[Going up works]
Let $Y\to X$ be a proper birational surjective morphism of normal proper varieties over a field $K$ of characteristic zero. If $X$ has the weak Hilbert property, then $Y$ has the weak Hilbert property. 
\end{lemma}
\begin{proof}
Let $Z\to Y$ be a ramified cover and consider the Stein factorization $Z'\to X$ of $Z\to Y\to X$. Note that the finite surjective morphism $Z'\to X$ is ramified.  (Indeed,  let $Z''\to Y$ be the pullback of $Z'\to X$ along $Y\to X$. Assume that $Z'\to X$ is \'etale.  Then $Z''\to Y$ is \'etale.   Moreover,  the finite surjective morphism $Z\to Y$ factors over the finite \'etale morphism $Z''\to Y$. Since the degree of $Z\to Y$ equals the degree of $Z'\to X$, we see that $Z\to Z''$ is of degree one, hence an isomorphism.  We conclude that $Z\to Y$ is \'etale.) 
In particular, since $X(K)$ is not strongly thin, there is a dense set $\Sigma$ of points in $Y(K)$ such that, for every $y$ in $\Sigma$, the fiber $Z_y$ does not have a $K$-point.
\end{proof}

For the reader's convenience, let us show that the weak Hilbert property for the symmetric product $\Sym^m(C\times \mathbb{P}^1_K)$ is a priori stronger than the weak Hilbert property for one of its desingularizations. In fact, one can not in general ``descend'' the weak Hilbert property along proper birational maps.

\begin{remark}[Going down fails]\label{remark:enriques}
Let $Y$ be the normal irreducible projective surface in $\mathbb{P}^3_{\mathbb{Q}}$ defined by
\[
x_0 x_2^4 + x_1 x_3^4 = x_0^2 x_1^3 + x_0^3x_1^2.
\]
Let $X\to Y$ be the minimal model of $Y$.  As explained in \cite[\S 3.2]{CZHP}, we have that $X$ is an Enriques surface. In particular, $X$ has the potential weak Hilbert property \cite{GvirtzChenMezzedimi}. However, the normal projective surface $Y$ is geometrically simply connected, but does not have the potential Hilbert property (see \cite[Theorem~1.3 and Remark~3.5]{CZHP}). In particular, as $Y$ obviously has a potentially dense set of rational points, this shows that the smoothness assumption is necessary in Conjecture \ref{conj:whp}. 
\end{remark}

\bibliographystyle{alpha}
\bibliography{puncture_refs}{}

\newcommand{\etalchar}[1]{$^{#1}$}
\def\cprime{$'$} \def\cprime{$'$}
\begin{thebibliography}{CDJ{\etalchar{+}}22}

\bibitem[AA03]{ArapuraArchava}
D.~Arapura and S.~Archava.
\newblock Kodaira dimension of symmetric powers.
\newblock {\em Proc. Amer. Math. Soc.}, 131(5):1369--1372, 2003.

\bibitem[Bar]{Bartsch}
F.~Bartsch.
\newblock New examples of geometrically special varieties: K3 surfaces,
  {E}nriques surfaces, and algebraic groups.
\newblock {\em Preprint}.

\bibitem[BCJW]{BCJW}
F.~Bartsch, F.~Campana, A.~Javanpeykar, and O.~Wittenberg.
\newblock The {W}eakly {S}pecial {C}onjecture contradicts orbifold {M}ordell,
  and thus abc.
\newblock {\em arXiv:2410.06643}.

\bibitem[BJR]{BJR}
F.~Bartsch, A.~Javanpeykar, and E.~Rousseau.
\newblock Weakly-special threefolds and non-density of rational points.
\newblock {\em Journal of the LMS, to appear. arXiv:2310.09065}.

\bibitem[BK05]{BrionKumar}
M.~Brion and S.~Kumar.
\newblock {\em Frobenius splitting methods in geometry and representation
  theory}, volume 231 of {\em Progress in Mathematics}.
\newblock Birkh\"auser Boston, Inc., Boston, MA, 2005.

\bibitem[Bog78]{BogomolovHol}
F.~A. Bogomolov.
\newblock Holomorphic tensors and vector bundles on projective manifolds.
\newblock {\em Izv. Akad. Nauk SSSR Ser. Mat.}, 42(6):1227--1287, 1439, 1978.

\bibitem[BSFP]{BSFPextensions}
L.~Bary-Soroker, A.~Fehm, and S.~Petersen.
\newblock Hilbert properties under base change in small extensions.
\newblock {\em Annali della Scuola Normale Superiore di Pisa, Classe di
  Scienze, to appear. arXiv:2312.16219}.

\bibitem[BSFP14]{BSFP}
L.~Bary-Soroker, A.~Fehm, and S.~Petersen.
\newblock On varieties of {H}ilbert type.
\newblock {\em Ann. Inst. Fourier (Grenoble)}, 64(5):1893--1901, 2014.

\bibitem[BT04]{BT}
F.~Bogomolov and Y.~Tschinkel.
\newblock Special elliptic fibrations.
\newblock In {\em The {F}ano {C}onference}, pages 223--234. Univ. Torino,
  Turin, 2004.

\bibitem[Cam]{CampanaBook-2}
F.~Campana.
\newblock Arithmetic {A}spects of {O}rbifold {P}airs.
\newblock {\em Chapter 2 in CRM Short Courses Springer \emph{Arithmetic
  geometry of logarithmic pairs and hyperbolicity of moduli spaces}.}

\bibitem[Cam04]{Ca04}
F.~Campana.
\newblock Orbifolds, special varieties and classification theory.
\newblock {\em Ann. Inst. Fourier (Grenoble)}, 54(3):499--630, 2004.

\bibitem[Cam11]{Ca11}
F.~Campana.
\newblock Orbifoldes g{\'e}om{\'e}triques sp{\'e}ciales et classification
  bim{\'e}romorphe des vari{\'e}t{\'e}s k{\"a}hl{\'e}riennes compactes.
\newblock {\em Journal de l'Institut de Math{\'e}matiques de Jussieu},
  10(4):809--934, 2011.

\bibitem[CCR22]{CCR}
B.~Cadorel, F.~Campana, and E.~Rousseau.
\newblock Hyperbolicity and specialness of symmetric powers.
\newblock {\em J. \'{E}c. polytech. Math.}, 9:381--430, 2022.

\bibitem[CDJ{\etalchar{+}}22]{CDJLZ}
P.~Corvaja, J.~L. Demeio, A.~Javanpeykar, D.~Lombardo, and U.~Zannier.
\newblock On the distribution of rational points on ramified covers of abelian
  varieties.
\newblock {\em Compos. Math.}, 158(11):2109--2155, 2022.

\bibitem[CDY]{CadorelDengYamanoi}
B.~Cadorel, Y.~Deng, and K.~Yamanoi.
\newblock Hyperbolicity and fundamental groups of complex quasi-projective
  varieties.
\newblock {\em arXiv:2212.12225}.

\bibitem[CHO76]{CampbellHowardOchiai}
L.~A. Campbell, A.~Howard, and T.~Ochiai.
\newblock Moving holomorphic disks off analytic subsets.
\newblock {\em Proc. Amer. Math. Soc.}, 60:106--108 (1977), 1976.

\bibitem[CO75]{CampbellOgawa}
L.~A. Campbell and R.~H. Ogawa.
\newblock On preserving the {K}obayashi pseudodistance.
\newblock {\em Nagoya Math. J.}, 57:37--47, 1975.

\bibitem[CW23]{CampanaWinkelmann}
F.~Campana and J.~Winkelmann.
\newblock Dense entire curves in rationally connected manifolds.
\newblock {\em Algebr. Geom.}, 10(5):521--553, 2023.
\newblock With an appendix by J\'{a}nos Koll\'{a}r.

\bibitem[CZ17]{CZHP}
P.~Corvaja and U.~Zannier.
\newblock On the {H}ilbert property and the fundamental group of algebraic
  varieties.
\newblock {\em Math. Z.}, 286(1-2):579--602, 2017.

\bibitem[CZ24]{ChenZhang}
S.~Chen and Z.~Zhang.
\newblock Weak approximation of symmetric products and norm varieties.
\newblock {\em Bulletin of the London Mathematical Society}, 56(5):1734--1747,
  2024.

\bibitem[Deb01]{DebarreBook}
O.~Debarre.
\newblock {\em Higher-dimensional algebraic geometry}.
\newblock Universitext. Springer-Verlag, New York, 2001.

\bibitem[Dem77]{Demazure}
M.~Demazure.
\newblock Automorphismes et d{\'e}formations des vari{\'e}t{\'e}s de {B}orel.
\newblock {\em Invent. Math.}, 39(2):179--186, 1977.

\bibitem[DG70]{SGA3}
M.~Demazure and A.~Grothendieck.
\newblock {\em Sch{\'e}mas en groupes I, II, III ({SGA} 3).}
\newblock Lecture Notes in Math. 151, 152, 153. Springer-Verlag, New York,
  1970.

\bibitem[EV92]{EVbook}
H.~Esnault and E.~Viehweg.
\newblock {\em Lectures on vanishing theorems}, volume~20 of {\em DMV Seminar}.
\newblock Birkh{\"a}user Verlag, Basel, 1992.

\bibitem[Fal83]{Faltings2-2}
G.~Faltings.
\newblock Endlichkeitss{\"a}tze f{\"u}r abelsche {V}ariet{\"a}ten {\"u}ber
  {Z}ahlk{\"o}rpern.
\newblock {\em Invent. Math.}, 73(3):349--366, 1983.

\bibitem[Fal94]{FaltingsLang}
G.~Faltings.
\newblock The general case of {S}. {L}ang's conjecture.
\newblock In {\em Barsotti {S}ymposium in {A}lgebraic {G}eometry ({A}bano
  {T}erme, 1991)}, volume~15 of {\em Perspect. Math.}, pages 175--182. Academic
  Press, San Diego, CA, 1994.

\bibitem[FJ74]{FreyJarden}
G.~Frey and M.~Jarden.
\newblock Approximation theory and the rank of abelian varieties over large
  algebraic fields.
\newblock {\em Proc. London Math. Soc. (3)}, 28:112--128, 1974.

\bibitem[GCM23]{GvirtzChenMezzedimi}
D.~Gvirtz-Chen and G.~Mezzedimi.
\newblock A {H}ilbert irreducibility theorem for {E}nriques surfaces.
\newblock {\em Trans. Amer. Math. Soc.}, 376(6):3867--3890, 2023.

\bibitem[GFP]{GFP}
N.~Garcia-Fritz and H.~Pasten.
\newblock Algebroid maps and hyperbolicity of symmetric powers.
\newblock {\em arXiv:2406.00835}.

\bibitem[Har77]{Har}
R.~Hartshorne.
\newblock {\em Algebraic geometry}.
\newblock Springer-Verlag, New York, 1977.
\newblock Graduate Texts in Mathematics, No. 52.

\bibitem[Has03]{HassettSurvey}
B.~Hassett.
\newblock Potential density of rational points on algebraic varieties.
\newblock In {\em Higher dimensional varieties and rational points ({B}udapest,
  2001)}, volume~12 of {\em Bolyai Soc. Math. Stud.}, pages 223--282. Springer,
  Berlin, 2003.

\bibitem[HH09]{smallness}
S.~Harada and T.~Hiranouchi.
\newblock Smallness of fundamental groups for arithmetic schemes.
\newblock {\em J. Number Theory}, 129(11):2702--2712, 2009.

\bibitem[HT00]{HassettTschinkel}
B.~Hassett and Y.~Tschinkel.
\newblock Abelian fibrations and rational points on symmetric products.
\newblock {\em Internat. J. Math.}, 11(9):1163--1176, 2000.

\bibitem[HT01]{HT}
B.~Hassett and Y.~Tschinkel.
\newblock Density of integral points on algebraic varieties.
\newblock In {\em Rational points on algebraic varieties}, volume 199 of {\em
  Progr. Math.}, pages 169--197. Birkh\"{a}user, Basel, 2001.

\bibitem[Iit82]{Iitaka}
S.~Iitaka.
\newblock {\em Algebraic geometry}, volume~76 of {\em Graduate Texts in
  Mathematics}.
\newblock Springer-Verlag, New York-Berlin, 1982.
\newblock An introduction to birational geometry of algebraic varieties,
  North-Holland Mathematical Library, 24.

\bibitem[Jav21]{JAut}
A.~Javanpeykar.
\newblock Arithmetic hyperbolicity: automorphisms and persistence.
\newblock {\em Math. Ann.}, 381(1-2):439--457, 2021.

\bibitem[Jav24]{Jnef}
A.~Javanpeykar.
\newblock Hilbert irreducibility for varieties with a nef tangent bundle.
\newblock {\em J. Ramanujan Math. Soc.}, 39(1):33--38, 2024.

\bibitem[JL24]{JLitt}
A.~Javanpeykar and D.~Litt.
\newblock Integral points on algebraic subvarieties of period domains: from
  number fields to finitely generated fields.
\newblock {\em Manuscripta Math.}, 173(1-2):23--44, 2024.

\bibitem[JR22]{JR}
A.~Javanpeykar and E.~Rousseau.
\newblock Albanese maps and fundamental groups of varieties with many rational
  points over function fields.
\newblock {\em Int. Math. Res. Not. IMRN}, (24):19354--19398, 2022.

\bibitem[KL]{KamenovaLehn}
L.~Kamenova and C.~Lehn.
\newblock Non-hyperbolicity of holomorphic symplectic varieties.
\newblock {\em Epiga, to appear. arXiv:2212.11411}.

\bibitem[Kob98]{KobayashiBook}
S.~Kobayashi.
\newblock {\em Hyperbolic complex spaces}, volume 318 of {\em Grundlehren der
  Mathematischen Wissenschaften [Fundamental Principles of Mathematical
  Sciences]}.
\newblock Springer-Verlag, Berlin, 1998.

\bibitem[KT02]{KrTsc}
A.~Kresch and Y.~Tschinkel.
\newblock Integral points on punctured abelian surfaces.
\newblock In {\em Algorithmic number theory ({S}ydney, 2002)}, volume 2369 of
  {\em Lecture Notes in Comput. Sci.}, pages 198--204. Springer, Berlin, 2002.

\bibitem[Luga]{Luger4}
C.~Luger.
\newblock Hilbert irreducibility for integral points on punctured linear
  algebraic groups.
\newblock {\em arXiv:2410.13403}.

\bibitem[Lugb]{Luger3}
C.~Luger.
\newblock A mixed fibration theorem for {H}ilbert irreducibility on non-proper
  varieties.
\newblock {\em Archiv der Mathematik, to appear. arXiv:2410.13741}.

\bibitem[Lugc]{Luger2}
C.~Luger.
\newblock Products of varieties with many integral points.
\newblock {\em International Journal of Number Theory, to appear.
  arXiv:2401.05203}.

\bibitem[Mat68]{Mattuck}
A.~Mattuck.
\newblock The field of multisymmetric functions.
\newblock {\em Proc. Amer. Math. Soc.}, 19:764--765, 1968.

\bibitem[MR22]{McKinnonRoth}
D.~McKinnon and M.~Roth.
\newblock Codimension two integral points on some rationally connected
  threefolds are potentially dense.
\newblock {\em J. Algebraic Geom.}, 31(2):345--386, 2022.

\bibitem[MR23]{MorrowRosso}
J.~S. Morrow and G.~Rosso.
\newblock A non-{A}rchimedean analogue of {C}ampana's notion of specialness.
\newblock {\em Algebr. Geom.}, 10(3):262--297, 2023.

\bibitem[MZ]{McKinnonZhu}
D.~McKinnon and Y.~Zhu.
\newblock The arithmetic puncturing problem and integral points.
\newblock {\em arXiv:1806.03180}.

\bibitem[PRT22]{PRT}
J.~V. Pereira, E.~Rousseau, and F.~Touzet.
\newblock Numerically non-special varieties.
\newblock {\em Compos. Math.}, 158(6):1428--1447, 2022.

\bibitem[RTW21]{RTW}
E.~Rousseau, A.~Turchet, and J.~T.-Y. Wang.
\newblock Nonspecial varieties and generalised {L}ang-{V}ojta conjectures.
\newblock {\em Forum Math. Sigma}, 9:Paper No. e11, 29, 2021.

\bibitem[Sch63]{Schwarzenberger}
R.~L.~E. Schwarzenberger.
\newblock Jacobians and symmetric products.
\newblock {\em Illinois J. Math.}, 7:257--268, 1963.

\bibitem[Ser92]{SerreTopicsGalois}
J.-P. Serre.
\newblock {\em Topics in {G}alois theory}, volume~1 of {\em Research Notes in
  Mathematics}.
\newblock Jones and Bartlett Publishers, Boston, MA, 1992.
\newblock Lecture notes prepared by Henri Damon [Henri Darmon], With a foreword
  by Darmon and the author.

\bibitem[Sik22]{Siksek}
S.~Siksek.
\newblock Integral points on punctured abelian varieties.
\newblock {\em Eur. J. Math.}, 8(suppl. 2):S687--S703, 2022.

\bibitem[Voj15]{VojtaExc}
P.~Vojta.
\newblock A {L}ang exceptional set for integral points.
\newblock In {\em Geometry and analysis on manifolds}, volume 308 of {\em
  Progr. Math.}, pages 177--207. Birkh\"{a}user/Springer, Cham, 2015.

\bibitem[Wit18]{Wittenberg18}
O.~Wittenberg.
\newblock Rational points and zero-cycles on rationally connected varieties
  over number fields.
\newblock In {\em Algebraic geometry: {S}alt {L}ake {C}ity 2015}, volume~97 of
  {\em Proc. Sympos. Pure Math.}, pages 597--635. Amer. Math. Soc., Providence,
  RI, 2018.

\bibitem[Wu]{Wu}
X.~Wu.
\newblock On a vanishing theorem due to {B}ogomolov.
\newblock {\em arXiv:2011.13751}.

\bibitem[Yam15]{Yamanoi1}
K.~Yamanoi.
\newblock Holomorphic curves in algebraic varieties of maximal {A}lbanese
  dimension.
\newblock {\em Internat. J. Math.}, 26(6):1541006, 45, 2015.

\end{thebibliography}
\end{document}